\documentclass[11pt,a4paper,twoside]{article}
\usepackage{amsthm,amssymb,amsfonts,amsmath}

\usepackage{xcolor,ulem}
\usepackage{enumitem}

\newtheorem{theorem}{Theorem}
\newtheorem{corollary}{Corollary}
\newtheorem{definition}{Definition}
\newtheorem{example}{Example}
\newtheorem{proposition}{Proposition}
\newtheorem{remark}{Remark}

\def\Ric{\operatorname{Ric}}
\newcommand\oc{\overset{\centerdot}}

\author{
 Sourav Nayak\footnote{Department of Mathematics, Indian Institute of Technology - Hyderabad, Sangareddy-502285, India
 \newline e-mail: {\tt ma22resch11004@iith.ac.in}.
 Orcid: 0009-0003-4330-8283}\ ,
 Dhriti Sundar Patra\footnote{Department of Mathematics, Indian Institute of Technology - Hyderabad, Sangareddy-502285, India
 \newline e-mail: {\tt dhriti@math.iith.ac.in, \tt dhritimath@gmail.com}.
 Orcid: 0000-0002-7958-1690}
 \ and
 Vladimir Rovenski \footnote{Department of Mathematics, University of Haifa, Mount Carmel, 3498838 Haifa, Israel
 \newline e-mail: {\tt vrovenski@univ.haifa.ac.il}.
 Orcid: 0000-0003-0591-8307.
 ({corresponding author})
 }
 }

\title{Characterizations of ${\mathcal S}$-Manifolds 
and \\ 
Splitting of Weak Almost ${\cal S}$-Manifolds}

\begin{document}

\date{}

\maketitle

\begin{abstract}
Weak metric structures, introduced by Rovenski and Wolak in 2022, extend Yano's $f$-structure and almost contact metric structure.
In this paper, we investigate curvature phenomena of weak almost ${\cal S}$-manifolds (w.a.$\,\cal S$-manifolds) focusing on the $f$-$(\kappa,\mu)$-nullity condition and its special case $R_{X,Y}\,\xi=0$. 
We establish several results that generalize known  rigidity theorems for almost ${\cal S}$-manifolds.
First, using the partial Ricci flow, we obtain dynamical characterizations of $\cal S$-manifolds:
starting from a w.a.$\,\cal S$-structure satisfying the curvature condition of $\cal S$-manifolds or the 
$f$-$(1,\mu)$-nullity condition, the flow evolves the structure exponentially fast toward an 
$\cal S$-structure. This extends results of Cappelletti Montano and Di~Terlizzi to the weak metric setting.
Next, we identify conditions under which a w.a.\,${\cal S}$-manifold admits a bi-Legendrian structure with totally geodesic foliations.
Finally, for w.a.$\,\cal S$-manifolds with $\kappa=\mu=0$, we prove a splitting theorem in which one factor is flat and the another has positive constant sectional curvature, generalizing 
classical results for almost $\cal S$-geometry.
These findings have consequences for the theory of Sasakian and $\cal S$- manifolds, the geometry of bi-Legendrian structures, and the behavior of weak metric contact manifolds under curvature constraints.



\end{abstract}



\section{Introduction}\label{sec:01}

Contact geometry 
plays an important role in many areas of mathematics and physics.
Within this framework, Yano's $f$-structure \cite{yan} on a differentiable manifold $M^{2n+s}$ -- defined using the conditions $f^3 + f = 0$ and rank$(f)=2n$  -- provides a unifying generalization of almost contact ($s=1$) and almost complex ($s=0$) structures. 
It arises naturally in the study of submanifolds of almost complex manifolds.
The~tangent bundle of an $f$-manifold splits into complementary subbundles: ${\cal D}:=f(TM)$ and $\ker f$, and when $\ker f$ is paralelizable, one obtains
 $f^2=-I+\sum\nolimits_{\,i}\eta^i\otimes\xi_i$,
where 
$\xi_1,\ldots,\xi_s$ are linear independent vector fields,
$\eta^1,\ldots,\eta^s$ are dual 1-forms,
and $I$ is the identity operator.
In this case, 
there exists an associated Riemannian metric $g$,
that is, 
\[
 g(f X, f  Y) = g(X, Y) - \sum\nolimits_{\,i}\eta^i(X)\eta^i(Y)
 \qquad(X,Y\in\mathfrak{X}_M),
\] 
$f$ is skew-symmetric,
and we obtain a metric $f$-structure 
$(f,\xi_i,\eta^i,g)$. 
Almost ${\cal S}$-manifolds (or, $f$-contact manifolds)
form a class of metric $f$-manifolds, satisfying 
 ${F}=d{\eta^1}=\ldots=d{\eta^s}$,
 where
 ${F}(X,Y):=g(X,{f} Y)$.
In this case, $\ker f$ is tangent to a $\mathfrak{g}$-foliation with flat totally geodesic leaves.
Almost ${\cal S}$-manifolds appear naturally in the study of symplectization, CR-integrability, the topology and dynamics of contact foliations (which generalize the Reeb 
flow on contact manifolds to higher dimensions), 
conjectures on the existence of closed leaves, $\mathfrak{g}$-foliations, Killing vector fields,  connections with skew-torsion symmetry, see \cite{AM-1995,BC-26,DIP-2001,Fin-2024,Fitz-2011,Goertsches-2}.
Their Lorentzian analogues also arise in models of spacetimes with electromagnetic field and globally hyperbolic geometry, see \cite{duggal}.

\smallskip

Several refinements of almost ${\cal S}$-manifolds have been introduced to capture curvature phenomena, see~\cite{b1970,CFF-1990,DTL-2006,DIP-2001,Fitz-2011,Goertsches-2}.
An~$f$-{K}-contact manifold is an almost ${\cal S}$-manifold, whose characteristic vector fields $\xi_1,\ldots,\xi_s$ are Killing (i.e., ${\cal L}_{\xi_i}\,g=0$, where ${\cal L}$ is the Lie~derivative).
Normal almost ${\cal S}$-manifolds -- called ${\cal S}$-manifolds -- generalize Sasakian manifolds; but unlike the Sasakian case, an $\cal S$-structure cannot exist on a simply connected compact manifold. 
Thus, the best structure we can hope to obtain on some manifolds is an almost $\cal S$-structure, see \cite{CMdiT-2007}.
Curvature conditions play a decisive role in characterizing these structures. For instance, ${\cal S}$-manifolds have constant $\xi$-sectional~curvature,
 $K(\xi,X)
 =1$, where 
 $\xi\in\ker f,\ X\perp\ker f$,
and are characterized among almost ${\cal S}$-manifolds by 
\begin{align}
\label{E-1k-0-nullity}
 R_{XY}\,\xi=\overline\eta(X) f^2Y 
 -\overline\eta(Y) f^2X\quad
 (X,Y\in\mathfrak{X}_M,\  \xi\in\ker f,\ \|\xi\|=1), 
\end{align} 
where 
$R_{X, Y}=[\nabla_X,\nabla_Y] -\nabla_{[X,Y]}$ is the curvature tensor, $\nabla$ is the Levi-Civita connection, 
and $\bar\eta=\sum_{\,i}\eta^i$.
Di~Terlizzi \cite{DTL-2006} showed that an almost $\mathcal{S}$-manifold 
satisfying the curvature-related restriction
\begin{align}\label{E-RXY-xi}
 R_{X,Y}\,\xi = 0 \qquad (X,Y\in\mathfrak{X}_M,\ \xi\in\ker f) 
\end{align}
splits locally 
as a product of 
Euclidean space 
and an $n$-sphere of radius $\frac{1}{2\sqrt{s}}$. 
Cappelletti Montano and Di~Terlizzi 
\cite{CM-2005,CMdiT-2007} later introduced the broader
$f$-$(\kappa,\mu)$-condition for the Reeb vector fields $\xi_i$,
\begin{equation}\label{E-k-mu-nullity}
 R_{X,Y}\,\xi_i =\kappa\{\overline\eta(X)f^2Y 
 -\overline\eta(Y)f^2X\} 
 +\mu\{\overline\eta(Y)h_i X -\overline\eta(X)h_i Y\}
 \quad (X,Y\in\mathfrak{X}_M,\ i=1,\ldots,s) ,
\end{equation}
where $\kappa,\mu\in\mathbb{R}$, and 
$h_i=\tfrac12\pounds_{\xi_i}{f}$, see \eqref{4.1},
which governs the behavior of 
bi‑Legendrian foliations.

\smallskip

A recent development by Rovenski and Wolak \cite{RWo-2} is the notion of weak metric structures, obtained by replacing the linear complex structure on the contact distribution with a non‑singular skew‑symmetric tensor, see Definition~\ref{Def-01}. 
Weak almost ${\cal S}$-manifolds 
(w.a.\,${\cal S}$-manifolds), see \cite{rst-43,rov-127}, generalize almost ${\cal S}$-manifolds and reveal new geometric phenomena.
%
For example, weak $f$-K-contact manifolds always have positive $\xi$-sectional curvature (for $f$-K-contact manifolds, this curvature is equal to 1), and a weak $f$-K-contact manifold equipped with a generalized Ricci soliton structure under certain conditions is an Einstein~manifold.
%
The $\,\cal S$-structure is rigid in the following sense:
a~w.a.$\,\cal S$-manifold is a weak $\,\cal S$-manifold if and only if it is an $\,\cal S$-manifold.

Advances in mathematics and physics
(such as
non-commutative geometry and superspace, string theory, 
the asymmetric energy-momentum tensor, and dark energy/matter)
have made the nonsymmetric metric tensor $G=g+F$, where $g$ is a pseudo-Riemannian metric (associated with gravity), and $F\ne0$ is a skew-symmetric tensor (associated with electromagnetism) more attractive than ever.
Weak metric structures naturally interact with nonsymmetric 
tensors $G=g+F$, relevant in mathematical physics,
and provide new tools for studying non-symmetric metric-affine manifolds, see~\cite{RZ-1}.


These advances raise 
natural questions:
\textit{Which w.a.$\,\mathcal{S}$-manifolds satisfyig the $f$-$(\kappa,\mu)$-nullity pro\-perty \eqref{E-k-mu-nullity} admit an $\,\mathcal{S}$-structure}? 
\textit{Under what conditions do they carry conjugate Legendrian foliations}?
\textit{When does a w.a.$\,\mathcal{S}$-manifold with curvature restriction \eqref{E-RXY-xi}~split}?

\smallskip

A key tool in our approach is the partial Ricci flow (PRF), introduced for foliated manifolds, see \cite[Section~5.3]{Rov-Wa-2021}. 
The PRF prescribes mixed sectional curvature 
and has recently been used to obtain dynamical characterizations of $f$-K-contact structures (K-contact structures when $s=1$), see \cite{rov-127} and discussion in Section~\ref{sec:02.2}.
Our main results 
include:

$\bullet$ \textbf{Dynamical characterization of ${\cal S}$-manifolds}: Starting from a w.a.$\,\mathcal{S}$-structure satisfying the curvature condition of $\,\mathcal{S}$-manifolds, the PRF evolves the structure exponentially fast toward an $\,\mathcal{S}$-structure,
this complements the results of the third author.

$\bullet$ \textbf{Bi‑Legendrian structures}:
We identify conditions under which w.a.\,${\cal S}$-manifolds admit conjugate Legendrian foliations, generalizing results of Cappelletti Montano and Di Terlizzi.

$\bullet$
\textbf{Scaling correspondence} between almost $\cal S$-structures and w.a.\,$\cal S$-structures satisfying 
\eqref{E-k-mu-nullity}.

$\bullet$ \textbf{Splitting theorems}: 
For w.a.\,$\cal S$-manifolds with $\kappa=\mu=0$, we prove a splitting theorem in which one factor is flat, generalizing classical results of Di Terlizzi.

These results show that weak metric structures preserve many of the essential 
features of 
contact and $f$-contact geometry while providing a broader and more flexible framework. They~also highlight the PRF as a powerful tool for understanding the evolution and rigidity of geometric structures.

\smallskip

 The paper is organized as follows. Section~\ref{sec:01} provides an introduction to the topic.
Section~\ref{sec:02} discusses known results 
and examples on w.a.\,${\cal S}$-manifolds and shows 
(in Proposition~\ref{L-cond-equiv}) that for a weak $f$-K-contact manifold satisfying  \eqref{E-nS-10}, the condition \eqref{E-1k-0-nullity} is invariant under the PRF.
This invariance is crucial for the dynamical characterization of 
${\cal S}$-manifolds. 
Sections~\ref{sec:main} and \ref{sec:04} present the main results and their proofs, 
see 
Theorems~\ref{P-K-S}--\ref{T-RXY=0}, their
Corollaries~\ref{C-K-S}--\ref{Cor-05} for $s=1$,
and three illustrative examples.
Section~\ref{sec:03} presents auxiliary statements 
concerning curvature and structure tensors of w.a.\,${\cal S}$-manifolds that generalize some results
of \cite{DTL-2006,tlk-2003,DIP-2001}.
Section~\ref{sec:concl} contains the conclusion.

\section{Preliminaries}
\label{sec:02}

Here we summarize the foundational notions used throughout the paper. Section~\ref{sec:02.1} reviews w.a.\,${\cal S}$-manifolds,
proves Proposition~\ref{L-cond-equiv} and presents two examples,
while Section~\ref{sec:02.2} recalls the partial Ricci flow (PRF) and its behavior on weak $f$-K-contact manifolds and proves Proposition~\ref{L-PRF-2}.


\subsection{Weak Almost $\,{\cal S}$-Manifolds}
\label{sec:02.1}

We extend the notion of Yano's $f$-structure
 (see \cite{yan}) and its special cases
 by replacing the linear complex structure on the contact distribution with a non‑singular skew‑symmetric tensor.

\begin{definition}\label{Def-01}
\rm
A \textit{weak metric $f$-structure}
on a differentiable manifold $M^{2n+s}$ consists of $(f,Q,\xi_i,\eta^i,g)$, where
$g$ is a Riemannian metric,
$f$ is a skew-symmetric (1,1)-tensor of rank $2\,n$, 
$\xi_i$ are mutually orthogonal unit vector fields, $\eta^i$ are dual 1-forms, i.e., 
$\eta^i(\xi_j)=\delta^i_j$, on $M$ satisfying
\begin{align}\label{2.2}
 g({f} X,{f} Y)&= g(X,Q\,Y) -\sum\nolimits_{\,i}{\eta^i}(X)\,{\eta^i}(Y)\qquad (X,Y\in\mathfrak{X}_M), 
\end{align}
where the self-adjoint tensor $Q:=-{f}^2 
+\sum\nolimits_{\,i}{\eta}^i\otimes {\xi}_i$ 
measures the deviation from classical metric
$f$-structure (where $Q=I$).
Hence, $f^3+f\,Q=0$.
\textit{Normality} of a~weak metric $f$-structure
is defined via the vanishing of the Nijenhuis torsion
 $N_f(X,Y) =
 {f}^2 [X,Y] + [{f} X, {f} Y] - {f}[{f} X,Y] - {f}[X,{f} Y]$
of~$f$. 
A~\textit{weak almost $\,{\cal S}$-structure} 
(briefly, a \textit{w.a.$\,{\cal S}$-structure}) 
is a weak metric $f$-structure satisfying 
\begin{align}\label{2.2-F}
 {F}=d{\eta^1}=\ldots=d{\eta^s},
 \quad 
 {\rm where}\quad {F}(X,Y):=g(X,{f} Y).
\end{align} 
When all vector fields ${\xi_i}$ are Killing,
the w.a.$\,{\cal S}$-structure is said to be \textit{weak $f$-K-contact}. 
\end{definition}

For $s=1$, weak metric $f$-manifolds, weak $f$-K-contact manifolds, and w.a.$\,\mathcal{S}$-manifolds are called weak metric manifolds, weak K-contact manifolds, and weak contact metric manifolds,~respectively.

\begin{remark}\rm
A weak metric $f$-structure is a~special case of a Riemannian almost product structure, defined by complementary orthogonal 
distributions ${\cal D}=f(TM)$ and $\ker f$, with Naveira's 36 distinguished classes.
Foliations appear when one or both distributions are involutive.
A~distribution 
is said to be {totally geodesic} if its second fundamental form va\-nishes, 
this is the case when {any geodesic of $M$ that is tangent to the distribution
at one point is tangent to the distribution at all its points}, e.g., \cite[Section~1.3.1]{Rov-Wa-2021}.
An involutive and totally geodesic distribution defines a totally geodesic~foliation. 
\end{remark}

We extend the concept of Legendrian distribution/foliation (see \cite{CM-2005,CMdiT-2007})
to weak metric $f$-manifolds.

\begin{definition}\label{D-002}\rm
An $n$-dimensional distribution $D$ on a w.a.$\,{\cal S}$-manifold $(M^{2n+s},f,Q,\xi_i,\eta^i,g)$ is called a \textit{Legendrian distribution} if $F(X,Y)=0$ for all vectors $X,Y$ in $D$. 
Two Legendrian distributions $D$ and $D'$
are said to be \textit{conjugate} if $f(D')=D$;
in this case, the tangent bundle splits as the orthogonal sum
$TM=D\oplus D'\oplus\ker f$. Conjugate distributions  form a pair known as a \textit{bi-Legendrian structure} on a w.a.$\,{\cal S}$-manifold.
When the Legendrian distribution is involutive, it defines a \textit{Legendrian foliation}.
\end{definition}

Putting $Y={\xi_i}$ in \eqref{2.2}, 
we get ${\eta^i}(X)=g(X,{\xi_i})$.
Thus, 
$f(TM)=\bigcap_{\,i}\ker{\eta^i}$ is true. 
By \eqref{2.2}, we get 
\begin{align*}
 {f}\,{\xi_i}=0,\quad {\eta^i}\circ{f}=0,\quad 
 Q\,\xi_i=\xi_i,\quad
 \eta^i\circ Q=\eta^i,\quad [Q,\,{f}]=0,\quad
 [\widetilde{Q},{f}]=0,\quad
 \eta^i\circ\widetilde{Q}=0 ,
\end{align*}
where
$\widetilde Q := Q-I$.
For a w.a.$\,\cal S$-manifold, 
$d{F} 
=0$ holds; thus, using the formula
 $3\,d{F}(X,Y,Z) =  X\,{F}(Y,Z) + Y\,{F}(Z,X) 
 + Z\,{F}(X,Y) - {F}([X,Y],Z) - {F}([Z,X],Y) 
 - {F}([Y,Z],X)$,
yields
\begin{align}\label{dPhi}
 (\nabla_X {F})(Y,Z)+(\nabla_Y {F})(Z,X)+(\nabla_Z {F})(X,Y)= 0 \quad (X,Y,Z \in \mathfrak{X}_M) .
\end{align}

The calculation of $\nabla f$ in Proposition~\ref{thm6.2AA} yields the new tensor $\mathcal{N}^{\,(5)}$ that complements the series 
of tensors $\mathcal{N}^{\,(1)},\,\mathcal{N}^{\,(2)}_i, \,\mathcal{N}^{\,(3)}_i,\,\mathcal{N}^{\,(4)}_{ij}$ well known in the classical theory $(Q=I)$, see~\cite{b1970}.

\begin{proposition}[see Corollary~1 in \cite{rst-43}]\label{thm6.2AA}
For a w.a.$\,{\cal S}$-structure, we have
\begin{align}\label{Eq-nabla-phi-N5}
 2\,g((\nabla_X f)Y,Z)=g(\mathcal{N}^{(1)}(Y,Z),f X) 
 + 2\,g(f X,f Y)\,\bar\eta(Z) 
 -2\,g(f X,f Z)\,\bar\eta(Y) +\mathcal{N}^{(5)}(X,Y,Z),  
\end{align}
where $\mathcal{N}^{(1)}:=N_f +2\sum\nolimits_{i}\,d\eta^i\otimes\xi_i$, and the $(0,3)$-tensor ${\cal N}^{\,(5)}$ is given by
\begin{align*}
\notag
 {\cal N}^{\,(5)}(X,Y,Z) &= f Z\,(g(X, \widetilde QY)) -f Y\,(g(X, \widetilde QZ)) +g([X, f Z], \widetilde QY)-\,g([X,f Y], \widetilde QZ)\\
 & +g([Y,f Z] -[Z, f Y] - f[Y,Z],\ \widetilde Q X) -X(g(\widetilde QY,fZ)).
\end{align*}
\end{proposition}


Several structural tensors appear in the weak setting. 
In particular, the tensors 
\begin{align}\label{4.1}
 h_i=\tfrac{1}{2}\pounds_{{\xi_i}}{f},\qquad
 \widetilde h_i:=Q^{-1} h_i\qquad (1\le i\le s). 
\end{align}
generalize the classical 
tensor $h$ of contact metric geometry and will be used repeatedly in what follows.

\begin{proposition}[see Theorem~2 in \cite{rst-43}]\label{thm6.2A}
For a w.a.$\,\mathcal{S}$-structure, we have
$d{\eta^i}({f} X,Y) = d{\eta^i}({f} Y,X)$ and $d{\eta^j}({\xi_i}, \cdot)=0$; moreover, $h_i$
vanishes if and only if $\,\xi_i$ is a Killing vector field.
\end{proposition}

For w.a.$\,\mathcal{S}$-manifolds, 
the distribution ${\cal D}=f(TM)$ is not involutive, and it is bracket-generating. Indeed, from  
$d{\eta^i}({f} X,X) = g({f} X,{f} X)$ and 
$2\,d{\eta^i}({f} X,X) 
= fX(\eta^i(X)) + g([X, {f} X], {\xi_i})$ we deduce
\begin{align}\label{E-non-inv}
 g([X, {f} X], {\xi_i}) = 2\,g({f} X,{f} X)>0\quad (X\in{\cal D}\setminus\{0\}).
\end{align}
For w.a.$\,\mathcal{S}$-manifolds, 
the distribution $\ker f$ is tangent to a $\mathfrak{g}$-foliation, see \cite{AM-1995,RWo-2}, with an abelian Lie~algebra: $[\xi_i,\xi_j]=0$.
By $d{\eta^j}({\xi_i},\,\cdot)=0$ (see Proposition~\ref{thm6.2AA}) we~have 
 $g(\nabla_{X}\,\xi_i,\,\xi_j)=
 -g(X,\nabla_{\xi_i}\,\xi_j)$.
Symmetrizing the above equality and using 
$g(\xi_i,\, \xi_j)=\delta_{ij}$ yield
$\nabla_{\xi_i}\,\xi_j+\nabla_{\xi_j}\,\xi_i = 0$.
From this and $[\xi_i, \xi_j]=0$ it follows that 
w.a.$\,\mathcal{S}$-manifolds satisfy the conditions
\begin{align}\label{R-03}
 g(\nabla_{X}\,\xi_i,\,\xi_j)=0,\qquad
 \nabla_{\xi_i}\,\xi_j = 0\quad (1\le i,j\le s).
\end{align} 
By \eqref{R-03}, 
$\ker f$ is tangent to a foliation with flat 
totally geodesic leaves.
By Corollary~1 in \cite{rst-43}, we have
\begin{align}\label{3.1AA}
 2\,g((\nabla_{\xi_i}{f})X,Y) 
 = g([\xi_i, {f} Y], \widetilde QX) -g([\xi_i,{f} X], \widetilde QY) - \xi_i(g(\widetilde QY,fZ)) 
 \quad (1\le i\le s).
\end{align}
Using \eqref{R-03}
and $f\xi_j=0$, we obtain 
$h_i\,\xi_j=[\xi_i, f\xi_j] 
-f[\xi_i,\xi_j]=0$; thus, 
$\widetilde h_i\xi_j=0$. 
For $X\in{\cal D}$, using \eqref{R-03},
we have 
$0 = 2{F}(\xi_i,X) = 2\,d\eta^j(\xi_i,X) 
 = g(\nabla_X\,\xi_j, \xi_i);$
thus, $g(\nabla_X\,\xi_j, \xi_i)=0\ (X\in\mathfrak{X}_M)$.
We~find
\begin{align}\label{4.2}
 (\pounds_{\xi_i}{f})Y
 = (\nabla_{\xi_i}{f})Y - \nabla_{{f} Y}\,\xi_i + {f}\nabla_Y\,\xi_i.
\end{align}
Using \eqref{4.2}
and $g((\nabla_{\xi_i}{f})Y, \xi_j)=0$, see \eqref{3.1AA} with $Z=\xi_j$, we get $\eta^j\circ h_i=0$; thus, $\eta^j\circ \widetilde h_i=0$,
\[
 (\eta^j\circ h_i)Y = 
 g((\pounds_{\xi_i}{f})Y, \xi_j)
 = g((\nabla_{\xi_i}{f})Y, \xi_j)
 -g(\nabla_{{f} Y}\,\xi_i , \xi_j)
 +g({f}\nabla_Y\,\xi_i, \xi_j) 
 = 0.
\] 

Denote by $^\top$ the ${\cal D}$-component of a vector. 
Our analysis is carried out under the following assumptions, which are trivially satisfied by almost $\,{\cal S}$-manifolds $(Q={I})$:
\begin{align}\label{Eq-2assump}
 \pounds_{\xi}\,{Q}&=0\qquad (\xi\in\ker f), \\
\label{E-nS-10}
 ( (\nabla_X Q)Y )^\top&=0 \qquad 
 (X, Y \in{\cal D}).
\end{align}

\begin{example}\label{Ex-dim3}\rm
$(i)$~Let $(M^{2n+s},\hat{f},{\xi_i},{\eta^i},\hat{g})$ be an almost $\,\mathcal{S}$-manifold, and $\lambda>0$ a differentiable function  on $M$. 
Then, $({f},Q,{\xi_i},{\eta^i},g)$ is 
a w.a.$\,\mathcal{S}$-structure on $M$ with $\mathcal{N}^{\,(5)}=0$, where
\begin{align*}
 {f} = \lambda^{1/2} \hat{f}, \quad
 Q = \lambda\,I + (1-\lambda)\sum\nolimits_i \eta^i\otimes\xi_i,
 \quad
 g = \lambda^{-1/2}\,\hat g 
 +(1-\lambda^{-1/2})\sum\nolimits_i\eta^i\otimes\eta^i;
\end{align*}
moreover, if $\lambda=const$, then the conditions
\eqref{Eq-2assump} and \eqref{E-nS-10} are true.

\smallskip

$(ii)$ Let $(M^{2+s},f,Q,\xi_i,\eta^i,g)$ be a 
w.a.$\,\mathcal{S}$-manifold. Then, $Q|_{\cal D} =-f^2|_{\cal D} ={\lambda}\,I_{\cal D}$ for some differentiable function ${\lambda}>0$ on $M$.
If \eqref{E-nS-10} is true, then $X({\lambda})=0$ for all $X\in{\cal D}$, and by Chow's theorem~\cite{SChow}, since the plane field $\cal D$ is bracket-generating, see \eqref{E-non-inv}, we get ${\lambda}=const\ne1$.
\end{example}

For an almost ${\cal S}$-structure, each 
(1,1)-tensor $h_i$ is self-adjoint, trace-free and anti-commutes with $f$, i.e., $h_i{f}+{f}\, h_i=0$, see \cite{DIP-2001}.
We generalized this result for a w.a.$\,\mathcal{S}$-structure.

\begin{proposition}
[see Proposition~4 of \cite{rst-43} and
Theorem~1 of \cite{rov-127}]
\label{C-2cond}
For a w.a.$\,\mathcal{S}$-manifold we have
\begin{align*}
 \nabla\, \xi_i = -f - fQ^{-1}h^*_i\quad (i=1,\ldots,s).
\end{align*} 
If a w.a.$\,\mathcal{S}$-manifold satisfies 
the condition \eqref{Eq-2assump}, then
each $(1,1)$-tensor $h_i$ is self-adjoint, traceless, anti-commutes with $f$ and commutes with $Q$; in addition, 
$\nabla_{\xi_i}Q=\nabla_{\xi_i}f=0$ and $\operatorname{tr} h_i=0$ hold, and
\begin{equation}\label{E-xi-n}
 \nabla\, \xi_i = -f - f\,\widetilde h_i\quad (i=1,\ldots,s) .
\end{equation} 
A w.a.$\,\mathcal{S}$-manifold is weak $f$-K-contact if and only if
\begin{equation}\label{E-xi-nK}
 \nabla\, \xi_i = -f\quad (i=1,\ldots,s) .
\end{equation} 
\end{proposition}

For a w.a.$\,\mathcal{S}$-manifold satisfying \eqref{Eq-2assump} and \eqref{E-nS-10}, using Proposition~\ref{C-2cond} and $\nabla\,Q=\nabla\,\widetilde Q$, we have
\begin{align}\label{E-nS-10b}
 (\nabla_X Q)Y 
 \notag 
 & = \big[(\nabla_X \widetilde Q)Y^\top 
 +\sum\nolimits_i\eta^i(Y) (\nabla_X \widetilde Q)\,\xi_i\big]^\top 
 +\sum\nolimits_i g((\nabla_X \widetilde Q)Y,\xi_i)\xi_i \\
&= -\sum\nolimits_i\big\{\eta^i(Y)\widetilde Q\nabla_X \xi_i 
 + g(\widetilde Q\nabla_X \xi_i,Y) \,\xi_i\big\}\qquad
 (X,Y \in\mathfrak{X}_M),
\end{align}
from which, using \eqref{R-03}, we obtain $\nabla_{\xi_j} Q=0$.

\begin{proposition}\label{L-cond-equiv}
(i)~Let a w.a.$\,\mathcal{S}$-manifold satisfy \eqref{Eq-2assump}, then the condition \eqref{E-1k-0-nullity} is equivalent to
\begin{align}\label{E-nabla-phi-nullity}
 (\nabla_X f)Y =\bar\eta(Y)f^2X -g(f^2Y,X)\bar\xi,\quad
 {\rm where}\quad \bar\xi=\sum\nolimits_{\,i}\xi_i,\ \
 \bar\eta=\sum\nolimits_{\,i}\eta^i;
\end{align}
moreover, if any of these conditions holds, then the structure is weak $f$-K-contact.

(ii) For an $f$-K-contact manifold satisfying \eqref{E-nabla-phi-nullity}, 
the condition \eqref{Eq-2assump} is true, and we have
\begin{align}\label{E-nS-10K}
 (\nabla_X Q)Y = \bar\eta(Y)\widetilde Q fX
 + g(\widetilde Q fX, Y)\,\bar\xi. 
 \end{align}
\end{proposition}

\begin{proof}

 \eqref{E-1k-0-nullity} 
 $\Rightarrow$ \eqref{E-nabla-phi-nullity}.
 For $Y\in{\cal D}$ and $X=\xi_i$, \eqref{E-1k-0-nullity} reads $R_{\,\xi_i,Y}\,\xi_i = -QY$.
 By this, \eqref{Eq-2assump} and \eqref{E-37} of Proposition~\ref{prop4.1}, we get $h_i^2 Y=0\ (Y\in{\cal D})$. 
 Since the tensor $h_i$ is self-adjoint, 
 see Proposition~\ref{C-2cond},
 and $h_i \xi_j=0$ is true, we get $h_i=\widetilde h_i=0\ (1\le i\le s)$. 
 By Proposition~\ref{thm6.2A}, $\xi_i\ (1\le i\le s)$ are Killing vector fields. Thus, $(f,Q,\xi_i,\eta^i,g)$
 is a weak $f$-$K$-contact structure,~and 
 we have, see \cite[Lemma~3]{rov-127},
 \begin{align}\label{E-R-varphi}
  R_{\,{\xi_i},X}=\nabla_X f\quad  (1\le i\le s) .
 \end{align} 
 From 
 \eqref{E-1k-0-nullity}, 
 \eqref{E-R-varphi} and the first Bianchi identity, we obtain
  $(\nabla_Xf)Y - (\nabla_Yf)X 
  =\overline \eta(Y) f^2X-\overline\eta(X) f^2Y$.
  Since $d{F}=0$, the equation \eqref{dPhi} yields
 \[
  g((\nabla_Z f)Y, X) = 
  g((\nabla_Xf)Y - (\nabla_Yf)X, Z)
  =g(\overline \eta(Y) f^2X-\overline\eta(X) f^2Y, Z) ,
 \]
 which reduces~to \eqref{E-nabla-phi-nullity}.

 \smallskip

\eqref{E-nabla-phi-nullity} $\Rightarrow$ 
\eqref{E-1k-0-nullity}.   
Differentiating $f\,\xi_i=0$ and using \eqref{E-xi-n} from Proposition~\ref{C-2cond} (requires \eqref{Eq-2assump}), 
 \eqref{E-nabla-phi-nullity}, and the definition of $Q$, gives $h_i\equiv 0\ (1\le i\le s)$. By this and
 Proposition~\ref{thm6.2A}, each $\xi_i$ is Killing and $(f,Q,\xi_i,\eta^i,g)$ is a weak $f$-$K$-contact structure. Thus, 
 \eqref{E-R-varphi} is true. Applying this and \eqref{E-nabla-phi-nullity}, we have
  $R_{\,\xi_i,X}Y =(\nabla_Xf)Y
  = \overline\eta(Y)f^2X -g(f^2Y,X)\overline\xi$.
 From this, using the first Bianchi identity, we get~\eqref{E-1k-0-nullity}.

 
 \smallskip


$(ii)$ 
For a weak  $f$-K-contact manifold, $\nabla \xi_i =-f$ for all $i$, and let us assume it satisfies \eqref{E-nabla-phi-nullity}.
 Taking the covariant derivative of 
 $Q=-f^2 + \sum\nolimits_{\,i}{\eta}^i\otimes {\xi}_i$ and using \eqref{E-nabla-phi-nullity} along with
 \eqref{E-xi-nK},
 we~get, 
\begin{align*}
 (\nabla_X Q)Y &= -(\nabla_X f^2)Y + \sum\nolimits_i\,\big[g(\nabla_X \xi_i,Y)\xi+ \eta(Y)\nabla_X \xi_i\big] \\
 &  = -(\nabla_X f)(fY)- f(\nabla_X f)Y 
 -\sum\nolimits_i\,\big[g(fX,Y) \xi_i +\eta^i(Y)\,fX\big] \\
 & =\sum\nolimits_i\,\big[\eta^i(Y)Q f X+g(QY, fX)\xi_i  
 -g(fX, Y) \xi_i - \eta^i(Y)\,fX\big]\\
 & = \sum\nolimits_i\,\big[\eta^i(Y)\widetilde Q\,f X +g(\widetilde QY, fX)\xi_i\big].
 \end{align*}
 For $X =\xi_j$, we have $\nabla_{\xi_j} Q=0$. Moreover, 
 \[
 (\mathcal{L}_{\xi_j}Q)X = (\nabla_{\xi_j} Q)X 
 -\nabla_{QX} \xi_j +Q\nabla_X \xi_j = fQX - QfX=0.
 \]
Therefore, for an $f$-K-contact manifold satisfying  \eqref{E-nabla-phi-nullity}, both \eqref{Eq-2assump} and \eqref{E-nS-10K} hold.
\end{proof}

\begin{example}
\rm
We construct a $(2+s)$-dimensional 
w.a.$\,\mathcal{S}$-manifold satisfying the condition
\eqref{E-k-mu-nullity}.
Let $\mathfrak{g}$ be a 
Lie algebra with a basis $\{e_1,e_2,\xi_1, \dots \xi_s\}$, and the Lie brackets
\begin{align*}
 [e_1,e_2]= 2\lambda^{1/2}\bar \xi, \quad 
 [e_2,\xi_i]=\lambda^{1/2}b_1\,e_1, \quad
 [\xi_i,e_1]=\lambda^{1/2}b_2\,e_2, \quad 
 [\xi_i, \xi_j]=0,
\end{align*}
where $b_1\ne b_2$ are real,
$\bar\xi=\sum_i \xi_i$,
and 
$\lambda>0$ is a real number. 
Define a weak metric $\,f$-structure $(f, Q,\xi_i,\eta^i,g)$ on the Lie group $G$, whose Lie algebra is $\mathfrak{g}$, as~follows:
\begin{align*}
 &f e_1=\lambda^{1/2} e_2, \quad 
 f e_2=-\lambda^{1/2} e_1, \quad 
 f\,\xi_i=0, \quad 
 Qe_1=\lambda\,e_1, \quad 
 Qe_2=\lambda\,e_2, \quad 
 Q\,\xi_i=\xi_i,
\end{align*}
$g$ is the left-invariant Riemannian metric such that $\{e_1,e_2,\xi_1, \dots \xi_s\}$ is an orthonormal frame and $\eta^i(X)=g(\xi_i,X)$ is true.
Note that the conditions \eqref{Eq-2assump} and \eqref{E-nS-10} are true, and
$(f, Q,\xi_i,\eta^i,g)$ is a w.a.$\,\cal S$-structure on $G$.
We have 
$\widetilde h_ie_1= \tfrac12(b_2-b_1)\,e_1, \ \widetilde h_ie_2=-\tfrac12(b_2-b_1)\,e_2$.
The following relations~hold:
\begin{align*}
 &\nabla_{e_1} e_2 =\tfrac12{\lambda^{1/2}}(2-b_1+b_2)\,\bar\xi,\quad 
 \nabla_{e_1}\,\xi_i= \tfrac12{\lambda^{1/2}}(b_1-b_2-2)\,e_2, \quad 
 \nabla_{e_2}\,\xi_i =\tfrac12{\lambda^{1/2}}(2+b_1-b_2)\,e_1, \\
 & \nabla_{\xi_i} e_1 = \tfrac12{\lambda^{1/2}}(b_1+b_2-2)e_2, \quad
 \nabla_{\xi_i} e_2 = \tfrac12{\lambda^{1/2}}(2-b_1-b_2)e_1, \quad \nabla_{e_1} e_1 = \nabla_{e_2}e_2 =\nabla_{\xi_i}\,\xi_j = 0.
\end{align*}
From the above relations, we find the components of the curvature tensor,
\begin{align*}
& R_{\,e_1,\,e_2}\,\xi_i=0, \quad  
 R_{\,e_1,\,\xi_i}\,\xi_j=-\big\{1-\tfrac14\,(b_2-b_1)^2\big\} f^2 e_1 + (2-b_1-b_2)\,h_je_1, \\ 
&R_{\,e_2,\, \xi_i}\,\xi_j= -\big\{1-\tfrac14(b_2-b_1)^2\big\} f^2 e_2 + (2-b_1-b_2)\,h_je_2.
\end{align*}
Therefore, 
$(G,f, Q, \xi_i, \eta^i, g)$ satisfies \eqref{E-k-mu-nullity}, with
 $\kappa = 1-\tfrac14 (b_2-b_1)^2\le 1$ and $\mu=2-b_1-b_2$.
\end{example}

\subsection{The Partial Ricci Flow}
\label{sec:02.2}

Let $(M,g)$ be a Riemannian manifold equipped with complementary orthogonal distributions,
$({\cal D},\widetilde{\cal D})$.
The \textit{partial Ricci flow} (PRF) of metrics on 
$(M,g,{\cal D},\widetilde{\cal D})$ is defined as,
see \cite[Section~5.3]{Rov-Wa-2021},
\begin{equation}\label{E-GF-Rmix-Phi}
 \oc g_t = -2\Ric_{\cal D}(g_t) +2\,s\,g^\top_t ,
\end{equation}
where $g_0=g$, $g^\top(X,Y):=g(X^\top,Y^\top)$
and $\,\oc{\phantom.}\,$ denotes the $t$-derivative.
The symmetric (0,2)-tensor
\begin{align}\label{E-rF-a}
 \Ric_{\cal D}(X,Y):=\sum\nolimits_{i} 
 g({R}_{X^\top,\,{\xi}_i}\,{\xi}_i,\,Y^\top),
\end{align}
where $^\top$ is the $\cal D$-component of the vector,
and $\xi_1,\ldots,\xi_s$ a local orthonormal frame of $\widetilde{\cal D}$,
is called the \textit{partial Ricci curvature} of ${\cal D}$.
The~dual to \eqref{E-rF-a} self-adjoint 
(1,1)-tensor,
called the \textit{partial Ricci~tensor}, is given by
 $\Ric_{\cal D}^\sharp(X) =\sum\nolimits_{\,i} \big(
 R_{X^\top,\,\xi_i}\,\xi_i\big)^\top$.
Since $X^\top=0\ (X\in{\cal \widetilde D})$, we have $\Ric_{\cal D}^\sharp(X)=0\ (X\in{\cal \widetilde D})$.
On~weak $f$-K-contact manifolds, the PRF simplifies dramatically:
the tensor $\Ric_{\cal D}^\sharp$
is positive definite on ${\cal D}=f(TM)$, 
using the orthoprojector ${I}_{\cal D}: TM\to{\cal D}$,
the flow reduces to ODE, see \cite{rov-127},
in~fact,
\begin{align}\label{E-PRF-fK}
 \oc\Ric\,\!_{\cal D}^\sharp = 4\Ric_{\cal D}^\sharp(\Ric_{\cal D}^\sharp -\,s\,{I}_{\cal D}) .
\end{align}
So the $f$-K-contact structure is a fixed point of~\eqref{E-GF-Rmix-Phi}: $\Ric_{\cal D}^\sharp=s\,{I}_{\cal D}$. We also have
\begin{align}\label{E-PRF-fT}
 \oc{T}\,^{\sharp}_\xi =
 2\,(\Ric_{\cal D}^\sharp - s\,{I}_{\cal D}) {T}^{\sharp}_\xi ,
\end{align}
where the 
(1,1)-tensor ${T}^{\sharp}_\xi:{\cal D}\to{\cal D}\ (\xi\in\ker f)$ is defined by
$g({T}^{\sharp}_\xi X,Y)=\tfrac12\,g([X,Y],\xi)$ for $X,Y\in{\cal D}$.
By the theory of ODEs, there exists a unique solution 
$\Ric_{{\cal D},t}^\sharp$ of \eqref{E-PRF-fK} for $t\in\mathbb{R}$;
thus, a solution $g_t$ of \eqref{E-GF-Rmix-Phi} exists for $t\in\mathbb{R}$ and it is unique.
%

\smallskip

The following theorem dynamically characterizes $f$-K-contact manifolds.
Namely, the PRF deforms a weak $f$-K-contact structure to a weak $f$-K-contact structure, whose limit is an $f$-K-contact structure.This deformation preserves all $\xi_i, \eta_i$, the distribution ${\cal D}$ and its orthogonality to $\ker f$.

\begin{theorem}[see \cite{rov-127}]\label{T-PRF-g2b}
Let $(M^{2n+s},f_0,Q_{\,0},\xi_i,\eta^i,g_0)$ be a weak $f$-K-contact manifold.
Then, there exist $\cal D$-adapted Riemannian metrics $g_t\ (t\in\mathbb{R})$
-- a solution of the PRF \eqref{E-GF-Rmix-Phi} -- such that each $(f_t,Q_{\,t},\xi_i,\eta^i,g_t)$ is a weak $f$-K-contact structure on $M$ with structural tensors defined on ${\cal D}$ as
\begin{equation}\label{E-Q-phi}
 Q_{\,t}|_{\cal D}=\tfrac1s\Ric_{{\cal D},t}^\sharp\,,
 \quad
 f_t\,|_{\,{\cal D}}=T_{\xi_i}^\sharp(t).
\end{equation}
Moreover,~$g_t$  converges exponentially fast, as $t\to-\infty$, to a limit metric $\hat g$
with $\widehat\Ric\,_{\cal D}^\sharp=s\,{I}_{\cal D}$ that provides an $f$-K-contact structure 
$(\hat f,\xi_i,\eta^i,\hat g)$ on $M$.
\end{theorem}

\begin{remark}\rm
By Theorem~\ref{T-PRF-g2b}, the results of Goertsches  and Loiudice \cite{Goertsches-2} on the topology of metric $f$-K-manifolds also hold for a wider class of weak $f$-K-contact manifolds.
\end{remark}

Using \eqref{E-Q-phi}, we have
\begin{align*}
 \oc g(X,fY) &= -2\Ric_{\cal D}(X,fY) +2\,s\,g^\top(X,fY) 
 =  -2\,s\,g(\widetilde QX,fY),\\
  g(X, \oc f Y) &
  = 2\,g(X, (\Ric\,\!_{\cal D}^\sharp -s\,I_{\cal D})fY)
  = 2\,s\,g(X, \widetilde QfY)
  = 2\,s\,g(\widetilde QX,fY).
\end{align*}
Thus, $(g(X,fY))'=\oc g(X,fY)+g(X, \oc f\,Y)=0$,
i.e., the PRF
preserves the equalities $d\eta^i=F$, see 
\eqref{2.2-F}.

\smallskip

A key observation for later use 
(in the proof of Theorem~\ref{P-K-S})
is that the curvature condition charac\-terizing 
${\cal S}$-manifolds,
is invariant under the PRF when the weak structure satisfies \eqref{Eq-2assump}.

\begin{proposition}\label{L-PRF-2}
For a weak $f$-K-contact manifold satisfying  \eqref{Eq-2assump},
the property \eqref{E-nabla-phi-nullity} remains invariant under the PRF \eqref{E-GF-Rmix-Phi}.
\end{proposition}

\begin{proof} 
In our case, $\oc g(X,Y) =-2\Ric_{\cal D}(X,Y) +2\,s\,g^\top(X,Y)$
and
$\oc\Ric\,\!_{\cal D}(X,Y) =
 \oc g(\Ric_{\cal D}^\sharp(X),Y) 
 +g(\oc\Ric\,_{\cal D}^\sharp(X),Y)$. 
Since $Q|_{\cal D}= {I}_{\cal D} +\widetilde Q$, where $\widetilde Q = Q-I$, by \eqref{E-Q-phi}, we have
\begin{equation}\label{E-nabla-Ric}
 \nabla_Y\Ric_{\cal D}^\sharp = 
 s\,\nabla_Y \widetilde Q + s\,\nabla_Y {I}_{\cal D}.
\end{equation}
Using \eqref{E-nabla-Ric}, $\nabla_Z\,\widetilde Q=\nabla_Z\,Q$ and \eqref{E-GF-Rmix-Phi}, we obtain
\begin{align*}
 (\nabla_Z\,\oc g)(X,Y) 
%
= -2\,g(\nabla_Z\,\Ric^\sharp_{\cal D}(X), Y)
 +2\,s\,g((\nabla_Z\,{I}_{\cal D})X, Y) = -2\,s\,g((\nabla_Z\,Q)X,Y).
\end{align*}
In view of \eqref{E-PRF-fT} and \eqref{E-Q-phi}, we have 
 $\oc\xi_{i} = \oc\eta\,^{i} = 0$ and
 $\oc f = \oc T\,\!_{\xi_i}^\sharp
 = 2\,(\Ric_{\cal D}^\sharp - s\,{I}_{\cal D}) {T}^{\sharp}_{\xi_i}$.
Using \eqref{E-PRF-fK}, \eqref{E-Q-phi} and 
$\oc\Ric\,\!_{\cal D}^\sharp =
\oc Q + \oc{I}_{\cal D}=\oc Q$,
we also have 
$ \oc{f^2} = -\oc Q
 = -\tfrac4s \Ric_{\cal D}^\sharp(\Ric_{\cal D}^\sharp 
 -\,s\,{I}_{\cal D}).$


Let's calculate the evolution under PRF of the equality \eqref{E-nabla-phi-nullity}. We have
\begin{align} \label{Eq-ev}
 (\oc\nabla_X f)Y +(\nabla_X \oc f)Y 
 =\bar\eta(Y)\oc{f^2}X
 -g(\oc{f^2}\,Y,X)\bar\xi
 -\oc g(f^2Y,X)\bar\xi .
\end{align}
By the above, $g((\nabla_X Q)Y,Z) = g((\nabla_X Q)Z,Y)$ (since $Q$ is self-adjoint), \eqref{E-nS-10K}
(see Proposition~~\ref{L-cond-equiv}),
and 
the general formula, e.g., \cite[Eq.~(4.8)]{Rov-Wa-2021},
\begin{equation*}
 2\,g(\oc\nabla_X Y, Z) = (\nabla_X\,\oc g)(Y,Z)
 +(\nabla_Y\,\oc g)(X,Z)-(\nabla_Z\,\oc g)(X,Y)
 \quad (X,Y,Z\in\mathfrak{X}_M),
\end{equation*}
we have
\begin{align*}
\notag
 g((\oc\nabla_X f)Y,Z)  &= 
 g(\oc\nabla_X(fY),Z) + g(\oc\nabla_X Y,fZ)\\
 & = \tfrac12\,\{(\nabla_X\,\oc g)(fY,Z)
 + (\nabla_{fY}\,\oc g)(X,Z) -(\nabla_Z\,\oc g)(X,fY)
 \notag\\
 & \quad +(\nabla_X\,\oc g)(Y,fZ) 
 +(\nabla_Y\,\oc g)(X,fZ) 
 -(\nabla_{fZ}\,\oc g)(X,Y)\} \notag\\ 
 & = -s\,\{g((\nabla_X Q)fY,Z) 
 +g((\nabla_{fY} Q)X,Z) 
 -g((\nabla_Z Q)X,fY) \notag \\
 &\quad 
 +g((\nabla_X Q)Y, fZ) + g((\nabla_Y Q)X,fZ)
 -g((\nabla_{fZ} Q)X,Y)\} \notag\\
 & =-2\,s\,\bar\eta(Y) g(\widetilde Q fX, fZ),  \\
 g( (\nabla_X \oc f)Y, Z)& = g(\nabla_X (\oc fY),Z) -g(\oc f\nabla_X Y,Z) \notag \\
 & = 2\,s\,g((\nabla_X Q)fY,Z)
 + 2\,g((\Ric_{\cal D}^\sharp 
 - s\,{I}_{\cal D}) (\nabla_X f)Y,Z) \notag \\ 
&  = 2\,s\,\bar\eta(Z) g(\widetilde QfX, fY) 
+ 2\,\bar \eta(Y) g((\Ric_{\cal D}^\sharp 
- s\,{I}_{\cal D})f^2X,Z) \notag \\
& = 2\,s\,\bar\eta(Z) g(\widetilde QfX, fY)
- 2\,s\,\bar\eta(Y) g(\widetilde Qf X,fZ). 
\end{align*}
Therefore, the inner product of the LHS of \eqref{Eq-ev} with $Z$ 
is as follows:
\begin{align}\label{E-PRF-cond1a}
 g( (\oc\nabla_X f)Y+(\nabla_X \oc f)Y, Z) =
 -4\,s\,\bar\eta(Y) g(\widetilde Q X, QZ)
 +2\,s\,\bar\eta(Z) g(\widetilde QX, QY).
\end{align}
Similarly, using $I_{\cal D}\,\widetilde Q=\widetilde Q$, \eqref{E-nS-10K} and \eqref{E-Q-phi}, for all $X,Y,Z\in\mathfrak{X}_M$ we calculate
\begin{align*}
\bar\eta(Y)g(\oc{f^2}\,X,Z) &=-\tfrac4s\,\bar\eta(Y) \, g(\Ric_{\cal D}^\sharp(\Ric_{\cal D}^\sharp 
 -\,s\,{I}_{\cal D})X,Z)
 = -4s\,\bar\eta(Y)g(Q^2 X-QX,Z) \\
 & =-4s\,\bar\eta(Y) g(\widetilde QX, QZ), \\
 -\bar\eta(Z)g(\oc {f^2}\,Y,X) &=\tfrac4s\,\bar\eta(Z) g(\Ric_{\cal D}^\sharp(\Ric_{\cal D}^\sharp 
 -s{I}_{\cal D})Y,X)=4\,s\,\bar\eta(Z)g( Q^2Y-QY,X) \\ 
 & = 4\,s\,\bar\eta(Z)g(\widetilde QY,QX), \\  
 -\bar\eta(Z)\oc g(f^2Y,X) &= 
 2\,\bar\eta(Z)\{\Ric_{\cal D}(f^2Y,X) 
 -s\,g^\top(f^2Y,X)\}  
 = 2\,s\,\bar\eta(Z)\{ g(Qf^2 Y,X) - g(f^2Y, X)\}  \\
 & = - 2\,s\,\bar\eta(Z) g(\widetilde QfX, fY) = - 2\,s\,\bar\eta(Z) g(\widetilde QX, QY) . 
\end{align*}
Therefore, the inner product of the RHS of \eqref{Eq-ev} with $Z\in\mathfrak{X}_M$ is as follows:
\begin{align}\label{E-PRF-cond1b}
g(\bar\eta(Y)\oc{f^2}X -g(\oc{f^2}\,Y,X)\bar\xi
 -\oc g(f^2Y,X)\bar\xi,\ Z) =
 -4\,s\,\bar\eta(Y) g(\widetilde Q X, QZ)
 +2\,s\,\bar\eta(Z) g(\widetilde QX, QY).
\end{align}
Since the RHS of \eqref{E-PRF-cond1a} and \eqref{E-PRF-cond1b}
coincide, \eqref{Eq-ev} is true. 
Hence, \eqref{E-nabla-phi-nullity} is preserved under the PRF.
\end{proof}

\section{Main Results}
\label{sec:main}

Our main results concern w.a.$\,\mathcal{S}$-manifolds $(M^{2+s},f,Q,\xi_i,\eta^i,g)$ satisfying the conditions 
\eqref{Eq-2assump} and~\eqref{E-nS-10}.


A family of Riemannian metrics $g_t$ on 
a w.a.$\,\mathcal{S}$-manifold $(M^{2n+s},f,Q,\xi_i,\eta^i,g)$
is said to be $\cal D$-\textit{adapted} if the distributions ${\cal D},\widetilde{\cal D}$ are $g_t$-orthogonal for all $t$ and
the metrics change on $\cal D$ only.


We obtain the dynamic characteristics of $\mathcal{S}$-manifolds among all w.a.$\,\mathcal{S}$-manifolds, complemen\-ting
\cite[Theorems~1.1 and~4.3]{CFF-1990}.
We show that the manifold admits an $\cal S$-structure obtained as the exponential limit under the PRF with ``$\cal D$-adapted" w.a.$\,\mathcal{S}$-structures satisfying conditions \eqref{E-nabla-phi-nullity} or~\eqref{E-1k-0-nullity}.

\begin{theorem}\label{P-K-S}
Let $(M^{2n+s},f_0,Q_0,\xi_i,\eta^i,g_0)$ be a w.a.$\,\mathcal{S}$-manifold satisfying the conditions \eqref{Eq-2assump} and~\eqref{E-1k-0-nullity}. Then, $(f_0,Q_0,\xi_i,\eta^i,g_0)$ is a weak $f$-$K$-contact structure with ${\cal N}^{\,(5)}=0$, and there exists a smooth family $(f_t,Q_{\,t},\xi_i,\eta^i,g_t)\ (t\in\mathbb{R})$
of weak $f$-K-contact structures on $M$ with 
$\cal D$-adapted metrics $g_t$ 
-- a solution of the PRF \eqref{E-GF-Rmix-Phi} --
that converges exponentially fast, as $t\to-\infty$, to an $\cal S$-structure.
\end{theorem}

By Theorem~\ref{P-K-S}, the known results on the topology of $\cal S$-manifolds also hold for a wider class of 
w.a.$\,\mathcal{S}$-manifolds satisfying the conditions \eqref{Eq-2assump} and~\eqref{E-1k-0-nullity}.

\begin{corollary}[see \cite{SDV-26}]\label{C-K-S}
Let a weak contact metric manifold $(M^{2n+1},f_0,Q_0,\xi,\eta,g_0)$ with 
the property 
\begin{align*}
 R_{XY}\,\xi&=\eta(X) f^2Y -\eta(Y) f^2X \quad(\xi\in\ker f,\ \|\xi\|=1), 
\end{align*} 
satisfy 
\eqref{Eq-2assump}. Then, $(f_0,Q_0,\xi,\eta,g_0)$ is a weak $K$-contact structure with ${\cal N}^{\,(5)}=0$, and there exists a smooth family $(f_t,Q_{\,t},\xi,\eta,g_t)\ (t\in\mathbb{R})$ of weak K-contact structures on $M$ with $\cal D$-adapted metrics $g_t$ 
-- a solution of the PRF \eqref{E-GF-Rmix-Phi} --
that converges exponentially fast, as $t\to-\infty$, to a Sasakian~structure.
\end{corollary}

For a w.a.$\,\mathcal{S}$-manifold satisfying $Q|_{\,\cal D} =\lambda\,I_{\,\cal D}$ for some 
real positive $\lambda\ne1$,
Theorem~\ref{P-K-S} is true without the requirement of \eqref{Eq-2assump}.
We illustrate this in the following.

\begin{example}
\rm
Let $(M^{2n+s},f_0,Q_0,\xi_i,\eta^i,g_0)$ be a 
w.a.$\,\mathcal{S}$-manifold
with $Q|_{\,\cal D} =\lambda\,I_{\,\cal D}$ for some 
real positive $\lambda\ne1$,
satisfy~the condition 
\eqref{E-1k-0-nullity}.  
From Theorem~\ref{T-PRF-g2b},
applying the theorem of uniqueness of a solution to ODE \eqref{E-PRF-fK}, we have a $\cal D$-conformal solution $\Ric^\sharp_{\mathcal D,t}
=\mu(t)\,I_{\cal D}$ of \eqref{E-Q-phi}. 
Since $\oc I_{\cal D}=0$, we have the ODE $\oc\mu=4\mu\,(\mu-s)$ with $\mu(0)=s\,\lambda>0$.
Thus, $\mu(t)=\frac{s}{1-C\,e^{4st}}\ (t\in\mathbb{R})$, where $C=1-\lambda^{-1}\ne0$.
Since our PRF \eqref{E-GF-Rmix-Phi} reduces to ODE; 
applying the theorem of uniqueness of a solution,
we~have $\oc g\,\!^\bot_t = 0$ and
$g\,\!^\top_t=\rho(t)\,g\,\!^\top_0$, where
$\rho(0)=1$ and
\begin{align*}
 &\oc g\,\!^\top_t =-2\mu(t)g^\top_t
 +2\,s\,g^\top_t = 2\,(s-\mu(t))\,g^\top_t
 \implies 
 \frac{d}{dt}\log \rho(t)=2\,(s-\mu(t)) 
 =  -\frac{2\,s\,Ce^{4st}}{1-C e^{4st}}, \\
 \implies & 
 \log\rho(t)= 
 -\int_0^t\frac{2\,s\,C e^{4s\tau}}{1-C e^{4s\tau}}\,d\tau
 =\frac{1}{2}\,\log\frac{1-C e^{4st}}{1-C} , 
 \quad \implies  
 \rho(t) = \Big(\frac{1-C e^{4st}} {1-C}\Big)^{1/2}.
\end{align*}
Thus, $g_t$ is conformal along ${\cal D}$ with the conformal factor equal to $\rho(t)$. 
By~the above, 
$f_t=f_0/\rho(t)$, $Q_{\,t}|_{\mathcal D}=\frac{1}{1-Ce^{4st}}\,I_{{\cal D}}$, and
$g_t=g^\bot_0+\rho(t)\,g^\top_0$.
Note that $\lim\limits_{t\to-\infty}Q_t=I_{\cal D}$.
By Theorem~\ref{P-K-S},
$(f_t,Q_{\,t},\xi_i,\eta^i,g_t)$ $(t\in\mathbb{R})$,
is a smooth family of weak $f$-K-contact structures on $M$ 
satisfying  \eqref{E-nabla-phi-nullity} or \eqref{E-1k-0-nullity}, 
that converges exponentially fast, as $t\to-\infty$, to an $\cal S$-structure
$(f_0/{\sqrt\lambda}, \xi_i,\eta^i, g^\bot_0+\sqrt{\lambda}\,g^\top_0)$.
\end{example}

The next result generalizes 
\cite[Lemma 1.1 and Propositions 1.1 and 1.2]{CMdiT-2007}.
Applying Theorem~\ref{P-K-S}, we show that a w.a.$\,\mathcal{S}$-manifold satisfying \eqref{E-k-mu-nullity} and \eqref{Eq-2assump} admits a ${\cal D}$-adapted $\mathcal{S}$-structure if $\kappa=1$, and has 
a bi-Legendrian structure
with totally geodesic foliations 
if $\kappa<1$ and \eqref{E-nS-10}~hold.

\begin{theorem}\label{T-04}
Let a w.a.$\,{\cal S}$-manifold $(M^{2n+s},f_0,Q_0,\xi_i,\eta^i,g_0)$ with the $f$-$(\kappa,\mu)$-nullity property 
\eqref{E-k-mu-nullity} satisfy the condition \eqref{Eq-2assump}. Then, $\kappa\le 1$,~and 

\smallskip 
(i) if $\kappa=1$, then 
$(f_0,Q_0,\xi_i,\eta^i,g_0)$ is a weak $f$-$K$-contact structure with ${\cal N}^{\,(5)}=0$,
and there exists a smooth family $(f_t,Q_{\,t},\xi_i,\eta^i,g_t)\ (t\in\mathbb{R})$
of weak $f$-K-contact structures on $M$ 
with $\cal D$-adapted metrics $g_t$
-- a solution of the PRF \eqref{E-GF-Rmix-Phi} --
that converges exponentially fast, as $t\to-\infty$, to an $\cal S$-structure.

\smallskip 
(ii) if $\kappa<1$ and the condition \eqref{E-nS-10} is true, then 
$\widetilde h_1=\ldots=\widetilde h_s$ and
{\rm Spec}$\,\widetilde h_i=\{0,\pm\sqrt{1-\kappa}\}$;
moreover, the orthogonal decomposition
$\,{\cal D}={\cal D}_+\oplus{\cal D}_-$ is true,
where the $n$-dimensional eigen-distributions ${\cal D}_\pm$ of the eigenvalues $\pm\sqrt{1-\kappa}$ determine conjugate Legendrian totally geodesic foliations.
\end{theorem}

\begin{corollary}[see \cite{SDV-26}]
Let a weak contact metric manifold $(M^{2n+1},f_0,Q_0,\xi,\eta,g_0)$ with the 
property 
\begin{equation}\label{E-k-mu-nullity-1}
 R_{X,Y}\,\xi =\kappa\{\eta(X)f^2Y -\eta(Y)f^2X\} 
 +\mu\{\eta(Y)h X -\eta(X)h Y\},
\end{equation}
where
$X,Y\in\mathfrak{X}_M,\ \xi\in\ker f$, and $\|\xi\|=1$,
satisfy the condition \eqref{Eq-2assump}. Then, $\kappa\le 1$,~and 

\smallskip 
(i) if $\kappa=1$, then 
$(f_0,Q_0,\xi,\eta,g_0)$ is a weak $K$-contact structure with ${\cal N}^{\,(5)}=0$,
and there exists a smooth family $(f_t,Q_{\,t},\xi,\eta,g_t)\ (t\in\mathbb{R})$
of weak K-contact structures on $M$ with 
$\cal D$-adapted metrics $g_t$
-- a solution of the PRF \eqref{E-GF-Rmix-Phi} --
that converges exponentially fast, as $t\to-\infty$, to a Sasakian structure.

\smallskip 
(ii) if $\kappa<1$ and the condition \eqref{E-nS-10} is true, then 
{\rm Spec}$\,\widetilde h=\{0,\pm\sqrt{1-\kappa}\}$;
moreover, 
the orthogonal decomposition
$\,{\cal D}={\cal D}_+\oplus{\cal D}_-$ is valid,
where the $n$-dimensional eigen-distributions ${\cal D}_\pm$ of the eigenvalues $\pm\sqrt{1-\kappa}$ determine conjugate Legendrian totally geodesic foliations (a bi-Legendrian structure). 
\end{corollary}


\begin{example}\label{C-K-S2}\rm
Let 
a w.a.$\,{\cal S}$-manifold $(M^{2n+s},f_0,Q_0,\xi_i,\eta^i,g_0)$ with $Q_0|_{\,\cal D} =\lambda\,I_{\,\cal D}$ for some 
real positive $\lambda\ne1$, satisfy 
\eqref{E-k-mu-nullity} and~\eqref{Eq-2assump}.
By Theorem~\ref{T-04}$(i)$,
we get $\kappa\le1$; moreover, 
for $\kappa=1$, we have $h_1=\ldots=h_s=0$. 
Thus, \eqref{E-k-mu-nullity} reduces to \eqref{E-1k-0-nullity}.
By Theorem~\ref{T-PRF-g2b}
the smooth family $(f_t,Q_{\,t},\xi_i,\eta^i,g_t)\ (t\in\mathbb{R})$,
where
\[
 f_t=\frac{f_0}{\rho(t)},\ \ 
 Q_{\,t}|_{\mathcal D}=\frac{1}{1-Ce^{4st}}\, I_{{\cal D}},\ \ 
 g_t=g^\bot_0+\rho(t)\,g^\top_0,\ \ 
 \rho(t) =\Big(\frac{1-C e^{4st}} {1-C}\Big)^{1/2},\ \
 C=1-\lambda^{-1}\ne0,
\]
consists of ${\cal D}$-adapted weak $f$-K-contact structures on $M$. 
By Theorem~\ref{T-04}(i), this family converges exponentially fast, as $t\to-\infty$, to an $\cal S$-structure
$(\lambda^{-1/2}f_0,\,\xi_i,\eta^i,\, g^\bot_0+\lambda^{1/2}\,g^\top_0)$ on $M$.
\end{example}

According to the following theorem, we can apply the transformation \eqref{E-homothety} to an almost $\,\cal S$-structure $(f',\xi_i,\eta^i,g')$ satisfying the curvature-related condition \eqref{E-RXY-xi}, and obtain examples of a w.a.$\,\cal S$-structure $(f,Q,\xi_i,\eta^i,g)$ satisfying the $(\kappa',\mu')$-condition \eqref{E-k-mu-nullity}
with $\kappa'=3\,{\lambda}^2(1-{\lambda})$ 
and $\mu'=2\,{\lambda}(1-{\lambda})$. 

\begin{theorem}\label{P-07}
Let $(M^{2n+s},f',\xi_i,\eta^i,g')$ be an almost $\,\cal S$-manifold, and $(f,Q,\xi_i,\eta^i,g)$ be a w.a.$\,\cal S$-structu\-re on $M$ such that, see Example~\ref{Ex-dim3}(i), 
\begin{align}\label{E-homothety}
 f={\lambda}\,f',\quad
 g={\lambda}^{-1} g'
 +(1-{\lambda}^{-1}) \sum\nolimits_i\eta^i\otimes\eta^i,\quad
 Q = \lambda^2 I + (1-\lambda^2)\sum\nolimits_i \eta^i\otimes\xi_i
\end{align}
for some real ${\lambda}>0$. Then, $(f,Q,\xi_i,\eta^i,g)$ satisfies the $(\kappa,\mu)$--nullity condition \eqref{E-k-mu-nullity}
if and only if $(f',\xi_i,\eta^i,g')$ satisfies the $(\kappa',\mu')$-condition 
with $\kappa'= {\lambda}^2(\kappa+3-{3\lambda})$ and $\mu'=\lambda(\mu+2-{2\lambda})$.
\end{theorem}

\begin{corollary}[see \cite{SDV-26}]
Let $(M^{2n+1},f',\xi,\eta,g')$ be an almost Sasakian manifold, and $(f,Q,\xi,\eta,g)$ be a weak almost contact metric structure on $M$ such that
\begin{align*}
 f={\lambda}\,f',\quad
 g={\lambda}^{-1} g'
 +(1-{\lambda}^{-1}) \eta\otimes\eta,\quad
 Q = \lambda^2 I + (1-\lambda^2)\eta\otimes\xi
\end{align*}
for some real ${\lambda}>0$. Then, $(f,Q,\xi,\eta,g)$ satisfies the $(\kappa,\mu)$--nullity condition \eqref{E-k-mu-nullity-1}
if and only if $(f',\xi,\eta,g')$ satisfies the $(\kappa',\mu')$-condition 
with $\kappa'= {\lambda}^2(\kappa+3-{3\lambda})$ and $\mu'=\lambda(\mu+2-{2\lambda})$.
\end{corollary}

The next result generalizes \cite[Theorem~2.1]{DTL-2006}.
We show that a w.a.$\,\mathcal{S}$-manifold 
with the property \eqref{E-RXY-xi} and some conditions for $Q$, splits, and $\mathbb{R}^{n+s}$ is one of its factors. 
In the $(2+s)$-dimensional case we find a 
condition for a w.a.$\,\mathcal{S}$-manifold to be flat, generalizing a result of Di~Terlizzi, see~\cite[Theorem~2.2]{DTL-2006}. 

\begin{theorem}\label{T-RXY=0}
Let a w.a.$\,\mathcal{S}$-manifold 
$(M^{2n+s},f,Q,\xi_i,\eta^i,g)\ (n>1)$ with the curvature-related property \eqref{E-RXY-xi} satisfy the conditions \eqref{Eq-2assump} and \eqref{E-nS-10}. 
Then, 
$\widetilde h_1=\ldots=\widetilde h_s$, 
and the manifold splits along the distributions 
\({\cal D}_-\oplus\ker f\) with flat leaves and \({\cal D}_+\), where 
the $n$-dimensional eigen-distributions ${\cal D}_\pm$ of the eigenvalues $\pm 1$, respectively, of $\widetilde h_i$ determine conjugate Legendrian totally geodesic foliations;
moreover, \(Q\,|_{\cal D}=\lambda I_{\cal D}\)\, 
for some constant
\(\lambda>0\), and integral manifolds of \({\cal D}_{+}\) have
constant sectional curvature $1/\sqrt{4\lambda}$.
If the curvature tensor of the leaves of ${\cal D}_+$ is zero at some point, then $n=1$.
\end{theorem}

\begin{example}
\rm
Let a w.a.$\,\mathcal{S}$-manifold $(M^{2+s},f,Q,\xi_i,\eta^i,g)$ with the curvature-related property \eqref{E-RXY-xi}, sa\-tisfy the condition \eqref{E-nS-10}. 
Then $Q|_{\mathcal D}=\lambda\,I_{\mathcal D}$ for some positive real
$\lambda\ne1$, see Example~\ref{Ex-dim3}(ii).
By~Theorem~\ref{T-RXY=0}, the conjugate Legendrian distributions ${\cal D}_+$ and ${\cal D}_-$ are one-dimensional. Fix some unit vector fields $e_1\in{\cal D}_+$ and $e_2\in{\cal D}_-$, 
hence, $fe_1=\lambda^{1/2} e_2$ and $f e_2=-\lambda^{1/2} e_1$. 
By~Theorem~\ref{T-RXY=0}, 
the manifold splits along \({\cal D}_-\oplus \ker f\) (with flat leaves) and ${\cal D}_{+}$, and the metric $g$ is flat. 
Then, we calculate
\begin{align*}
 &g(\nabla_{e_2} e_1 , e_2) =- g(\nabla_{e_2} e_2, e_1)=0, \quad 
 g(\nabla_{e_2} e_1,\xi_i) = -g(\nabla_{e_2} \xi_i, e_1)=0,  \\ 
 & g(\nabla_{\xi_i}\,e_1, e_2) =- g(\nabla_{\xi_i}\, e_2, e_1)=0,\quad g(\nabla_{\xi_i}\,e_1, \xi_j) 
 =- g(\nabla_{\xi_i}\,\xi_j, e_1)=0, \\ 
 & g(\nabla_{e_1} e_1, e_1)=0, \quad 
 g(\nabla_{e_2} e_1, e_1) = g(\nabla_{\xi_i}\,e_1, e_1)=0\qquad (1\le i,j\le s) .
\end{align*}
From the above and $\bar\xi=\sum_{\,i}\xi_i$, we have the relations,
\begin{align*}
 \nabla_{e_1}\,e_2=2\,\lambda\,\bar\xi, \quad 
 \nabla_{e_2}\,e_1=0, \quad 
 \nabla_{e_1}\,\xi_i = -2\,\lambda\,e_2, \quad 
 \nabla_{e_2}\,\xi_i = \nabla_{\xi_i} e_1 
 =\nabla_{\xi_i}\,\xi_j=0.
\end{align*}
Therefore, the subspaces ${\cal D}_- \oplus \operatorname{ker} f$ rotate along the ${\cal D}_+$-curves (geodesics) with the speed $2\lambda$.
\end{example}


\begin{corollary}[see \cite{SDV-26}]\label{Cor-05}
Let \(M^{2n+1}(\varphi,Q,\xi,\eta,g)\) be a weak contact metric manifold satisfying conditions \(\eqref{E-RXY-xi}\), 
\(\eqref{Eq-2assump}\), and \(\eqref{E-nS-10}\).
Then \(M\) splits 
along
\({\cal D}_{-}\oplus\operatorname{span}\{\xi\}\) and
\({\cal D}_{+}\), where the former distribution is flat and
\(({\cal D}_{+},{\cal D}_{-})\) defines a bi-Legendrian structure.
If \(n>1\), then \(Q\,|_{\cal D}=\lambda I_{\cal D}\)\, 
for some constant
\(\lambda>0\), and integral manifolds of \({\cal D}_{+}\) have
constant sectional curvature $1/\sqrt{4\lambda}$.
\end{corollary}

\begin{remark}\rm
By a transformation analogous to \eqref{E-homothety}, the w.a.$\,\mathcal{S}$-structure 
in Theorem~\ref{T-RXY=0} (and its corollary) reduces to the almost $\,\mathcal{S}$-structure. 
\end{remark}

\section{Proof of Main Results}
\label{sec:04}

In this section, we prove Theorems~\ref{P-K-S}--\ref{P-07}, using the results of Sections~\ref{sec:02} and~\ref{sec:03}.

\begin{proof}[\textbf{Proof of Theorem~\ref{P-K-S}}] 
By Proposition~\ref{L-cond-equiv}, the conditions 
\eqref{E-1k-0-nullity} and \eqref{E-nabla-phi-nullity} for a w.a.$\,\mathcal{S}$-manifold are equiva\-lent, and $(f_0,Q_0,\xi_i,\eta^i,g_0)$ is an $f$-$K$-structure.
Using the Levi-Civita connection $\nabla$ of $g$, we have 
\begin{align*}
 N_f(X,Y) = ({f}\nabla_Y{f} 
 -\nabla_{{f} Y}{f}) X 
 -({f}\nabla_X{f} 
 -\nabla_{{f} X}{f}) Y.
\end{align*}
From this
and \eqref{E-nabla-phi-nullity}, we~find 
\begin{align*}
 g(\mathcal{N}^{(1)}(Y,Z),f_0 X) 
 = g(f_0 (\nabla_Z f) Y -(\nabla_{f_0Z} f_0) Y
 -f_0 (\nabla_Y f_0)Z +(\nabla_{f_0Y} f_0) Z,\ f_0X) =0. 
\end{align*}
Applying this and \eqref{E-nabla-phi-nullity}, we reduce \eqref{Eq-nabla-phi-N5} of Proposition~\ref{thm6.2A} to the equality 
${\cal N}^{\,(5)}=0$.
By Theorem~\ref{T-PRF-g2b} and Proposition~\ref{L-PRF-2}, there exist a $\cal D$-adapted family $(f_t,Q_{\,t},\xi_i,\eta^i,g_t)$ of weak $f$-K-contact structures on $M$, which converges exponentially fast, as $t\to-\infty$, to a limit $f$-K-contact structure 
$(\hat f,\xi_i,\eta^i,\hat g)$ satisfying~\eqref{E-nabla-phi-nullity}.
According to the results of \cite{BC-26}, 
$(\hat f,\xi_i,\eta^i,\hat g)$ is an $\cal S$-structure.
\end{proof}

\begin{proof}[\textbf{Proof of Theorem~\ref{T-04}}] 

Using the assumptions and Proposition~\ref{C-2cond}, we derive
\begin{align*}
 R_{\xi_i,X}\,\xi_j \notag 
 & = \nabla_{\xi_i}((f+f\widetilde h_j)X) 
 + (f+f\widetilde h_j)(f+f\widetilde h_i)X \\  
 & = -fQ^{-1}(\nabla_{\xi_i}h_j)X +(f^2-f^2\widetilde h_j+f^2\widetilde h_i 
 -f^2\widetilde h_j \widetilde h_i)X.
\end{align*} 
Using this, we find
$Q R_{\xi_i,X}\xi_j - f R_{\xi_i, fX}\xi_j = 2 (Q f^2+Q^2 \widetilde h_j \widetilde h_i)X.$
From 
\eqref{E-k-mu-nullity} we obtain
$ Q R_{\xi_i,X}\xi_j - f R_{\xi_i, fX}\xi_j = 2\,\kappa\,Q f^2 X$.
Comparing both of these equations, we find
 $h_jh_i = (\kappa-1) Q f^2= h_ih_j$.   
Taking $i=j$, yields 
\[
 \widetilde{h}_i^2 |_{\mathcal D} = Q^{-2}h_i^2|_{\mathcal D}=(1-\kappa)I_{\mathcal D};
\] 
hence, $\kappa\le1$. Since $\widetilde h_i$ is self-adjoint, we have $\widetilde h_1=\ldots=\widetilde h_s=0$.

$(i)$ This reduces \eqref{E-xi-n} to $\nabla\,\xi_i=-f$, see \eqref{E-xi-nK}. 
By Proposition~\ref{thm6.2A}, 
$g(\nabla_X\,\xi_i,Y)+g(\nabla_Y\,\xi_i,X)=0$, 
i.e., $\xi_i$ is a Killing vector.
Hence, our manifold is weak $f$-K-contact. 
By \eqref{E-k-mu-nullity} with $\kappa=1$ and $h_i=0$, we get \eqref{E-1k-0-nullity} which is equivalent to \eqref{E-nabla-phi-nullity} by Proposition~\ref{L-cond-equiv}. 
By Theorem~\ref{P-K-S}, there exist a ${\cal D}$-adapted fami\-ly $(f_t,Q_{\,t},\xi_i,\eta^i,g_t)$ of weak $f$-K-contact structures on $M$,
which converges exponentially fast, as $t\to-\infty$, to a limit ${\cal S}$-structure
$(\hat f,\xi_i,\eta^i,\hat g)$.

$(ii)$ 
Let $\kappa<1$, then the eigenvalues of $\widetilde h_i$ are $\{0, \pm\sqrt{1-\kappa}\}$. Fix $i$ and $x\in M$.
By Proposition~\ref{prop-5.1}(ii), we have ${\cal D}_x={\cal D}^+_x\oplus{\cal D}^-_x$, where ${\cal D}^+_x$ (resp., ${\cal D}^+_x$) consists of the eigenvectors of $\widetilde h_i$ 
with positive (respectively, negative) eigenvalues. 
Any vector $X\in T_xM$ can be decomposed as $X=X^+ + X^-$, thus 
$\widetilde h_i X^\pm = \pm\sqrt{1-\kappa}\,X^\pm$. 
Using $X^+ - X^- =\frac1{\sqrt{1-\kappa}}\,\widetilde h_i X$, we~derive
\begin{align*}
 \widetilde h_j X =  \widetilde h_j(X^+ + X^-) 
 = \tfrac1{\sqrt{1-\kappa}}(\kappa-1) Q^{-1}f^2(X^+ - X^-) = Q^{-2}(-f^2) h_iX = \widetilde h_iX
\end{align*}
for all $j=1,\ldots,s$.
Therefore, 
$\widetilde h_1=\ldots=\widetilde h_s:=\widetilde h$.
Since the self-adjoint operator $\widetilde h$ has three distinct eigenvalues $\{0, \pm\sqrt{1-\kappa}\}$, the tangent bundle is decomposed into eigen-distributions: $\operatorname{ker} f$, ${\cal D}_+$, and ${\cal D}_-$, associated, respectively, with each eigenvalue. 
Using $f\widetilde h = -\widetilde hf$,
see Proposition~\ref{C-2cond},
we obtain $\widetilde{h}(fX)=-\sqrt{1-\kappa}\,fX\ (X\in{\cal D}^+)$, which ensures
$fX\in{\cal D}^-$. Thus, $f$ maps ${\cal D}^+$ onto ${\cal D}^-$ and vice versa, so $\dim {\cal D}^+ = \dim {\cal D}^- = n$. 
Consequently, the contact distribution decomposes as a direct sum
${\cal D}={\cal D}_+ \oplus{\cal D}_-$
of orthogonal $n$-dimensional eigen-distributions ${\cal D}_\pm$ corresponding to the eigenvalues $\pm\sqrt{1-\kappa}$ of~$\widetilde{h}_i$.
The distribution $\operatorname{ker} f$ is integrable, since $[\xi_i,\xi_j]=0$ for all $1\leq i\le j\le s$. 

Let $X,Y \in {\cal D}_+$ (${\cal D}_-$, respectively), then using \eqref{E-xi-n} we get 
$ g(\nabla_X \xi_i ,Y)= -(1\pm\sqrt{1-\kappa})\,g(f X,Y) =0=g(\nabla_Y \xi_i, X) \ (i=1,\dots ,s).$
Thus, $0=2\,g(X,f Y)= 2\,d\eta^i(X,Y)=- \eta^i([X,Y])$ for all $i=1,\dots,s$. 
It follows from \eqref{E-k-mu-nullity} that $R_{X,Y}\, \xi_i=0$
for $X,Y \in {\cal D}$. Using this for $X,Y,Z\in {\cal D}_+$ ($X,Y,Z\in{\cal D}_-$, respectively), we acquire
\begin{align*}
 0&= g(R_{X,Y}\,\xi_i,Z) 
 = -(1\pm\sqrt{1-\kappa})\,(\nabla_Z\,{F})(X,Y)\mp 2\sqrt{1-\kappa}\,g(f [X,Y],Z).
\end{align*}
Applying \eqref{Eq-phi02} and \eqref{Eq-phi001} 
of Proposition~\ref{prop-5.1} in the above relation yields, we get $g(f [X,Y],Z)=0$. 
Using this and $\eta^i([X,Y])=0$ for all $i=1,\dots,s$, we have $[X,Y]\in {\cal D}_+$ \ ($[X,Y]\in{\cal D}_-$, respectively),
which shows that the distributions ${\cal D}^\pm$ are involutive.

Next, fix $X,Z\in{\cal D}_-$ and $Y\in {\cal D}_+$, then using the fact $R_{e_1,e_2}\,\xi_i=0\ (e_1, e_2 \in{\cal D})$, we have 
\begin{align*}
 0&= g(R_{X,Y}\,\xi_i ,Z) = 
 (1+\sqrt{1-\kappa})\big\{ (\nabla_X {F})(Y,Z)
 + (\nabla_Y {F})(Z,X)\big\} 
 - 2\,\sqrt{1-\kappa}\, g(f \nabla_Y X,Z), \\
 & = (1+\sqrt{1-\kappa}) (\nabla_Z {F})(X,Y)- 2\,\sqrt{1-\kappa}\, g(f \nabla_Y X,Z).
\end{align*}
Applying \eqref{Eq-phi02}, gives 
$(\nabla_Z {F})(X,Y)=0$, hence $g(\nabla_Y X,f Z)=0$ for $X,Z\in{\cal D}_+$ and $Y\in{\cal D}_-$. 
Thus $\nabla_Y X\ 
(X\in{\cal D}_+,\,Y\in{\cal D}_-)$ 
is orthogonal to ${\cal D}_-$. 
Taking  \(Y,Z\in{\cal D}_+\) and \(X\in{\cal D}_-\), and using \eqref{E-xi-n}, gives
\[
g(\nabla_Y Z, \xi_i) = -g(\nabla_Y \xi_i, Z) = (1+ \sqrt{1-\kappa})\,g(f Y, Z) = 0,\quad
 g(\nabla_Y Z, X) = -g(Z, \nabla_Y X) = 0.
\]
By the above, 
\(\nabla_Y Z \perp {\cal D}_- \oplus  \operatorname{ker }f\); hence \(\nabla_Y Z\in{\cal D}_+\), i.e., \({\cal D}_+\) is totally geodesic. Applying the same procedure for $X,Z\in{\cal D}_+$ and $Y\in {\cal D}_-$ and using \eqref{Eq-phi02}
of Proposition~\ref{prop-5.1}, we conclude that ${\cal D}_-$ defines a totally geodesic foliation.
According to Definition~\ref{D-002}, ${\cal D}_\pm$
define conjugate Legendrian foliations.
\end{proof}

\begin{proof}[\textbf{Proof of Theorem~\ref{P-07}}] 
Note that $g'(X,\xi_i)=\eta^i(X)=g(X,\xi_i)$ and 
$g'(X,Y^\top)={\lambda}\,g(X,Y^\top)$, where $Y^\top$ is the ${\cal D}$-component of $Y$.
Let $\nabla'$ and $\nabla$ be the Levi-Civita connections of the metrics $g'$ and $g$, respectively. 
In the Koszul formula for $g'$,  we expand each term using 
$g'(X,Y) = {\lambda}\,g(X,Y) + (1-{\lambda})\,\sum\nolimits_i\eta^i(X)\eta^i(Y)$, see \eqref{E-homothety}, to obtain
\begin{align*}
 2\,g'(\nabla'_X Y, Z) 
 =&\, 2\,{\lambda}\,g(\nabla_X Y,Z) 
 + (1-{\lambda})
 \sum\nolimits_i\big\{2\eta^i(\nabla_X Y) \eta^i(Z) 
 +(\nabla_X \eta^i)(Y)\eta^i(Z) 
  +\eta^i(Y)(\nabla_X \eta^i)(Z) \\ 
 & 
 +(\nabla_Y\eta^i)(X)\eta^i(Z) +\eta^i(X)(\nabla_Y \eta^i)(Z) -(\nabla_Z \eta^i)(X)\eta^i(Y)-\eta^i(X)(\nabla_Z \eta^i)(Y)\big\}.
\end{align*}
Simplifying the above by using 
the relation $(\nabla_X\,\eta^i)(Y)= g(\nabla_X\,\xi_i,Y) =-g(f  X +f \widetilde h_i X,Y)$, see \eqref{E-xi-n}, and the equality 
$g'(\nabla_X Y,Z) = {\lambda}\,g(\nabla_X Y,Z) 
+(1-{\lambda})\sum\nolimits_i\eta^i(\nabla_X Y) \eta^i(Z) $, we have 
\begin{align*}
 &2\,g'(\nabla'_X Y, Z) = 2\,g'(\nabla_X Y, Z) 
 +(1-{\lambda})\sum\nolimits_i\big\{-\eta^i(Z) g(f  X 
 {+}f \widetilde h_i X,Y) -\eta^i(Y) g( f  X 
 {+}f \widetilde h_i X,Z) \\
 & {-}\eta^i(Z) g(f  Y {+}f  \widetilde h_i Y,X) 
 -\eta^i(X)g(f  Y {+}f  \widetilde h_i Y, Z) 
 +\eta^i(Y)g(f  Z {+}f  \widetilde h_i Z,X) +\eta^i(X)g(f  Z {+}f  \widetilde h_i Z,Y) \big\}.
\end{align*}
From this, using the 
equalities $h_i=h_j:=h\ (1 \le i,j\le s)$ for a weak $(\kappa,\mu)$-manifold, we obtain
\begin{align}\label{E-nabla'}
 \nabla'_X Y =\nabla_XY+ (1-{\lambda})\,\{ g(\widetilde hX, fY)\bar\xi-\bar\eta(Y) fX-\bar\eta(X) fY\}.
\end{align}
 From \eqref{E-k-mu-nullity} we conclude that
$ R_{\,\xi_i,X}Y = \kappa \{ g(QX,Y)\bar\xi 
- \bar\eta(Y)QX \} 
 + \mu \{ g(hX,Y)\bar\xi - \bar \eta(Y)hX \}$. Using the previous equation and the definition of $Q$, \eqref{E-05b} of Proposition~\ref{prop-5.2} gives
\begin{equation*}
(\nabla_{hX}f )Y 
= {\lambda}^4(1-\kappa)\big\{ g(X,Y)\bar\xi-\bar\eta(Y)X \big\} 
+ {\lambda}^2\big\{\bar\eta(Y)hX - g(hX,Y)\bar\xi\big\} + \lambda^2\sum\nolimits_i \eta^i(X)[\eta^i(Y)\bar\xi- \bar\eta(Y)\xi_i].
\end{equation*}
In the above equation, we substitute $hX$ in place of $X$ and use the equality 
 $h^2 =(1-\kappa)\big\{Q^2 
 -\sum\nolimits_i\eta^i\otimes\xi_i\big\}$,
see the proof of Theorem~\ref{T-04}(i), to obtain
\begin{equation}\label{Eq-varphi}
 (\nabla_X f )Y = \lambda^2[g(X,Y)\bar\xi-\bar\eta(Y) X +\sum\nolimits_i \eta^i(X)\{\bar \eta(Y) \xi_i-\eta^i(Y) \bar\xi\} +g(hX,Y)\bar\xi-\bar\eta(Y) hX.
 \end{equation}
 Set 
 $A(X,Y) := g(\widetilde hX,f  Y)\xi - \bar\eta(Y)f  X - \bar\eta(X)fY$,
and using this along with \eqref{E-nabla'}, we find the relation
\[
 R'_{X,Y}\,\xi_i = R_{X,Y}\,\xi_i 
 +(1-{\lambda})\big\{(\nabla_X A)(Y,\xi_i) 
 -(\nabla_Y A)(X,\xi_i) +A(X,A(Y,\xi_i)) -A(Y,A(X,\xi_i)) \big\}.
\]
We have $A(X, A(Y,\xi_i)) = -A(X, f Y)  =
 -g( hX, Y)\,\bar\xi+\bar\eta(X)f^2 Y$. Using this and  \eqref{Eq-varphi}, we have 
 \begin{align*}
 R'_{X,Y}\,\xi_i &= R_{X,Y\,}\xi_i + (1-{\lambda}) 
 \{3(\bar\eta(X) f ^2Y-\bar\eta(Y) f ^2X ) + 2(\bar\eta(Y) hX -\bar\eta(X) hY)\}.
\end{align*}
Using the definition of $f'$ and equalities
$h_i = \tfrac12\mathcal{L}_{\xi_i}\,f 
= \tfrac12\,{\lambda}\,\mathcal{L}_{\xi_i}\,f'
= {\lambda}\,h_i'$
and \eqref{E-k-mu-nullity}, we have
\begin{align*}
 R'_{X,Y}\,\xi_i=
 {\lambda}^2(\kappa+3-3{\lambda})\{\bar\eta(X) f'^2Y
 -\bar\eta(Y) f'^2X\} +\lambda(\mu+2-2{\lambda})
 \{\bar\eta(Y) h_i'X-\bar\eta(X) h_i'Y\}.
\end{align*}
Taking $\kappa'={\lambda}^2(\kappa+3-3{\lambda})$ and $\mu'={\lambda}(\mu+ 2-2{\lambda})$, we see that $(f ',\xi_i,\eta^i,g')$ is a $(\kappa',\mu')$-manifold. 
\end{proof}

\begin{proof}[\textbf{Proof of Theorem~\ref{T-RXY=0}}] 

By Theorem \ref{T-04}$(ii)$, 
$\widetilde h_1=\ldots=\widetilde h_s$, and
the distributions ${\cal D}_{\pm}$ are integrable and totally geodesic. 
According to Definition~\ref{D-002}, ${\cal D}_\pm$
determine conjugate Legendrian 
foliations.

We~have $\nabla_{[\xi_j,X]}\,\xi_i=-R_{\xi_j,X}\,\xi_i=0\ (X\in{\cal D}^-)$; hence, $\widetilde h_i[\xi_j,X]=-[\xi_j,X]$, that is, $[\xi_j,X]\in{\cal D}^-\ (1\le j\le s,\ X\in{\cal D}^-)$. By this and $[{\xi_i},\xi_j]=0$, see \eqref{R-03}, the distribution ${\cal D}^-\oplus\ker f$ is involutive. 
Therefore, we can consider a foliated chart $U$ with coordinates $x_1,\ldots,x_{2n+s}$ such that $\{\partial_j\}_{j>n}$ is a local basis of
${\cal D}^-\oplus\ker f$. There exist differentiable functions $c^j_a$ on $U$ such that 
$X_a=\partial_a+\sum_{j>n}c^j_a \partial_j$ 
is a basis of ${\cal D}^+$.
Since $[\partial_j, X_a]\in{\cal D}^-\oplus\ker f$ for
$j>n$ and $a\le n$, we can write 
$[\partial_j, X_a]=X+\sum_{j\le s}\sigma^j\xi_j$, where $X\in{\cal D}^-$ and $\sigma^j$ are differentiable functions on $M$. 
From 
 $\nabla_{[\partial_j, X_a]}\,\xi_i = \nabla_X\,\xi_i 
 +\sum\nolimits_{j\le s}\sigma^j\nabla_{\xi_j}\,\xi_i = 0$ 
we conclude that
 $\xi_i$ is parallel along $[\partial_j, X_a]$.
Then, using \eqref{E-xi-n} and the condition \eqref{E-RXY-xi}, we get
\begin{align*}
 0 = \nabla_{[\partial_j, X_b]}\,\xi_i 
 = \nabla_{\partial_j}\nabla_{X_b}\,\xi_i 
 -\nabla_{X_b}\nabla_{\partial_j}\,\xi_i
 = -2\,\nabla_{\partial_j}(f X_b),
\end{align*}
and, since $f X_a\in{\cal D}^-$, we obtain
 $\nabla_{f X_a}(f X_b)=0$.
From the above, we conclude that the integral manifolds of ${\cal D}^-\oplus\ker f$ are flat and totally geodesic. By the de Rham decomposition theorem, \(M\) is locally a Riemannian product, and one of its factors is locally isometric to 
\(\mathbb{R}^{n+s}\).

\smallskip

In the rest of the proof, all calculations are done for $X,Y,Z,{V}\in{\cal D}_+$.  Using \eqref{E-xi-n}, we have
\[
 g((\nabla_X\, f )Y,\xi_i)
 = -g((\nabla_X\, f  )\xi_i,Y)
 = 2\,g( f  X, f  Y)
 = 2\,g(QX,Y)\quad (i=1,\dots,s) .
\]
Hence, from \eqref{Eq-phi02} of Proposition~\ref{prop-5.1}, we acquire
\begin{align}\label{E-phi2}
 (\nabla_X\, f )Y=2\,g(QX,Y)\overline\xi.
\end{align}
As ${\cal D}_{+}$ defines a totally geodesic foliation, by \eqref{E-phi2} we have
\begin{align} \label{E-RXY1}
  g(\nabla_X \nabla_Y  f  Z, f  {V})-g(\nabla_X \nabla_Y Z, {V}) \notag &=g(\nabla_X (2g(QY,Z)\overline\xi
 + f \nabla_Y Z), f  {V}) -g(\nabla_X \nabla_Y Z, {V})\notag \\
 &
 = -4s\,g(QY,Z)g(QX,{V}) + g(\nabla_X \nabla_Y Z, \widetilde Q{V}).
\end{align}
Also, from \eqref{E-phi2}, we have 
 $g(\nabla_{[X,Y]}\, f  Z,  f  {V})
 -g(\nabla_{[X,Y]} Z,{V}) 
   = g(\nabla_{[X,Y]} Z,\widetilde Q{V}).$
Combining this with \eqref{E-RXY1}, we get
\begin{align*}
 g(R_{X,Y}  f  Z,  f  {V}) - g(R_{X,Y}\,Z, Q{V})  
 = 4s\big\{g(Y,Q{V})g(QX,Z) - g(QY,Z)g(X,Q{V})\big\}.
\end{align*}
Hence, $g(R_{X,Y}fZ, f{V})=0$ for $X,Y,Z,{V} \in {\cal D}_{+}$, indicating that  after replacing $QW$ by $W$, we have
\begin{align}\label{E-RXY3}
    g(R_{X,Y}\,Z,W) = 4\,\{g(QY,Z)\,g(X,W) -g(QX,Z)\,g(Y,W)\} \qquad
 (X,Y,Z,W\in{\cal D}_+).
\end{align}
Using the skew-symmetry of the Riemannian curvature tensor, we have 
\[
 0=g(R_{X,Y}\,Z,Z)= 4\,\{g(QY,Z)\,g(X,Z) -g(QX,Z)\,g(Y,Z)\}.
\]
Taking $Z\in{\cal D}_+$ as a unit vector and setting $X=Z$ yields  
\begin{align*}
    g(QY,Z)=g(QZ,Z)\,g(Y,Z)\qquad\text{for all }Y\in\mathcal D_+ .
\end{align*}
By symmetry of $Q$, we have 
$QZ = g(QZ,Z)\,Z$. 
Thus, every nonzero vector in ${\cal D}_+$ is an eigenvector of $Q$. 
Hence, $Q|_{{\cal D}_+} =\lambda\, I|_{{\cal D}_+}$ for some smooth function $\lambda$ on $M$. Since $\varphi$ and $Q$ commute and $\varphi ({\cal D}_-)={\cal D}_+$, for $V \in {\cal D}_ -$ we have
\[
\varphi QX =Q\varphi V = \lambda\, \varphi V \implies \varphi (QV -\lambda \, V)=0.
\]
Since, $\varphi$ is injective on ${\cal D}$, it follows that $QV =\lambda \, V$. From this and above, $Q|_{{\cal D}} =\lambda\, I|_{{\cal D}}$.
The compatibility conditions \eqref{Eq-2assump} and \eqref{E-nS-10} force $\lambda $ to be a positive real constant. Therefore, \eqref{E-RXY3} simplifies~as
\[
g(R_{X,Y}\,Z,W) = 4\lambda\,\{g(Y,Z)\,g(X,W) -g(X,Z)\,g(Y,W)\} \qquad
 (X,Y,Z,W\in{\cal D}_+).
\]
Hence, the integral manifolds of ${\cal D}_+$ have constant sectional curvature $1/\sqrt{4\lambda}$ and the manifold is locally isometric to $\mathbb{R}^{n+1}\times \mathbb{S}^n(4\lambda)$.


Suppose, 
that 
the curvature tensor of the leaves of ${\cal D}_+$ is zero at some point.
Since the manifold satisfies \eqref{E-RXY-xi},
\eqref{Eq-2assump} and \eqref{E-nS-10},
the equation~\eqref{E-RXY3} is valid
and reduces to 
\begin{align}\label{E-cor1}
 g(Y,Z)\,g(X,W) = g(Y,W)\,g(X,Z) \quad
 (X,Y,Z,W\in {\cal D}_+).
\end{align}
Take $Z=Y$ and $W=X$ in \eqref{E-cor1}, then  
$g(X,Y)^2= g(X,X)\,g(Y,Y)$.
By the Cauchy-Schwarz inequality, this equality holds if and only if $X$ and $Y$ are linearly dependent. 
Therefore,
$\dim{\cal D}_+=n=1$.
\end{proof}


\section{Auxiliary Results}
\label{sec:03}

Here, we study many aspects of the curvature and structure tensors of the w.a.$\,\mathcal{S}$-structure satisfying conditions \eqref{Eq-2assump}, \eqref{E-nS-10}, and generalize the corresponding results known in the classical case $Q=I$. 

The following result generalizes Proposition 3.1 and its corollary of \cite{tlk-2003}.

\begin{proposition}\label{prop4.1}
For a w.a.$\,\mathcal{S}$-manifold satisfying \eqref{Eq-2assump}, we have the following:
\begin{align}\label{E-37}
 QR_{\xi_i, X}\xi_j-f R_{\xi_i,\,f X} \xi_j=2(h_j h_iX+Qf^2 X) .
\end{align}
\end{proposition}

\begin{proof}
Let's compute $R_{\xi_i\, X}\xi_j $ by applying \eqref{E-xi-n} (from Proposition~\ref{thm6.2A}) along with $\nabla_{\xi_i}\,\xi_j=0$, see \eqref{R-03}:
\begin{align}\label{E-38}
 R_{\xi_i, X}\xi_j
 & = -\nabla_{\xi_i} ( f X +f\,Q^{-1}h_jX) 
 +f[\xi_i,X] + f\,Q^{-1}h_j[\xi_i,X]. 
\end{align}
Applying $f$ to both sides of \eqref{E-38} and then recalling Proposition~\ref{C-2cond},
we have
\begin{align*}
& f R_{\xi_i, X}\xi_j
  =-f^2\big[Q^{-1}(\nabla_{\xi_i}\,h_j)X 
  +\nabla_X \xi_i 
 + Q^{-1}h_j\nabla_X \xi_i\big] \notag \\
 & = (\nabla_{\xi_i}\,h_j)X +Q\nabla_X\,\xi_i +h_j\nabla_X \xi_i \sum\nolimits_{\,k=1}^s\big[\eta^k((\nabla_{\xi_i}\,h_j)X) +\eta^k(\nabla_X\,\xi_i)
 -\eta^k(Q^{-1}h_j\nabla_X\,\xi_i)\big]\xi_k. 
\end{align*}
In the above equation, applying \eqref{E-xi-n} and the fact that $\eta^k((\nabla_{\xi_i}\,h_j)X)=0$ (follows by taking the covariant derivative of $g(h_jX,\xi_i)=0$ along $\xi_i$), we achieve the relation 
\begin{align}\label{E-36} 
(\nabla_{\xi_i}\,h_j)X = f R_{\xi_i, X}\xi_j +h_iX
-h_jX + Qf X - f\,Q^{-1}h_jh_iX .
\end{align}
From \eqref{E-36}, we have the following:
\begin{align*}
 & QR_{\xi_i, X}\xi_j =-f (\nabla_{\xi_i}\,h_j)X -h_iX+h_jX+Qf^2X +h_jh_iX,\\
 & f R_{\xi_i, f X}\xi_j = -f (\nabla_{\xi_i}\,h_j)X -h_iX+h_jX-Qf^2X  -h_jh_iX.
\end{align*}
Combining these two relations, we acquire \eqref{E-37}. 
\end{proof}

The following result generalizes Proposition 2.5 of \cite{DIP-2001} and is used in Proposition~\ref{prop-5.2}.

\begin{proposition}\label{P-05}
For a w.a.$\,\mathcal{S}$-manifold satisfying \eqref{Eq-2assump} and \eqref{E-nS-10}, the 
following is true:
\begin{align}\label{E-04}
 (\nabla_X f) Y + (\nabla_{f X}f) f Y = &\,2\,g(fQX,fY) \overline{\xi}+\overline{\eta}(Y)f^2QX -\sum\nolimits_j \eta^j(Y)Qh_jX  \notag \\
& -\tfrac{1}{2}\big\{\widetilde Q(\nabla_X f) Y 
 +(\nabla_{\widetilde QX} f) Y \big\}
-\tfrac{1}{2}\sum\nolimits_j g((Q+h_j)Y,\widetilde QX)\xi_j .
\end{align}
\end{proposition}

\begin{proof}
 The covariant derivative of ${F}(Y,Z) =g( Y, f Z)$ along $X$ gives $(\nabla_X {F})(Y,Z)=g((\nabla_X f)Z,Y)$. Using this, \eqref{E-xi-n} and \eqref{E-nS-10b}, we have 
\begin{align}
 (\nabla_X {F})(f Y,Z) & - (\nabla_X {F})( Y,f Z) 
 = g((\nabla_Xf)Z ,f Y) -g((\nabla_X f^2)Z,Y)
 +g(f(\nabla_X f)Z,Y) \notag \\
 & = -\sum\nolimits_j\,\big[\eta^j(Y)g((Q+h_j)X,fZ) +\eta^j(Z)g((Q+h_j)X,fY)\big] \label{E-39}.
\end{align}
Replacing $Z$ by $f Z$ in \eqref{E-39}, we acquire that 
 \begin{align}
  (\nabla_X {F})(fY, fZ) - (\nabla_X {F})(Y,f^2 Z)
  = \sum\nolimits_j \big[\eta^j(Y)g((Q+h_j)X,QZ) 
  -\overline{\eta}(Y)\eta^j(X) \eta^j(Z)\big]. \label{E-40}
 \end{align}
 A simple computation gives  
\[
 (\nabla_X {F})( Y,f^2 Z) =-(\nabla_X {F})( Y, Q Z) -\sum\nolimits_j\big[\eta^j(Z)g((Q+h_j)X,Y)
 -\overline{\eta}(Z)\eta^j(X)\eta^j(Y)\big].
\]
Inserting this in \eqref{E-40}, gives
\begin{align}
 (\nabla_X {F})(f Y, f Z) &+ (\nabla_X {F})( Y, Z)\notag\ = \sum\nolimits_j\big[
 \eta^j(Y)g((Q+h_j)X,QZ)   
 -\eta^j(Z)g((Q+h_j)X,Y) \notag\\ 
 & +\, \overline{\eta}(Z)\eta^j(X)\eta^j(Y) 
  -\overline{\eta}(Y)\eta^j(X) \eta^j(Z)\big] 
  -(\nabla_X {F})( Y,\widetilde Q Z) . \label{E-41}
\end{align}
Now, since $d{F}=0$ (as ${F}=d\eta^i$ for $1\le i\le s$), using \eqref{dPhi},
where $(\nabla_X {F})(Y,Z)=g((\nabla_X f)Z,Y)$ for all $X,Y\in TM$, we have
\begin{align*}
 &\quad (\nabla_X {F})(Y,Z)+ (\nabla_Y {F})(Z,X)+(\nabla_Z {F})(X,Y) \\ 
 &+(\nabla_{f X} {F})(f Y, Z) +(\nabla_{f Y} {F})(Z,f X)  +(\nabla_Z {F})(f X, f Y)\\
 & +(\nabla_{f X} {F})(Y,f Z) 
 +(\nabla_Y {F})(f Z, f X)  
 +(\nabla_{f Z} {F})(f X,Y) \\
 &-(\nabla_X {F})(f Y, f Z)
 -(\nabla_{f Y} {F})(f Z, X)-(\nabla_{ f Z} {F})(X, f Y) =0.
\end{align*}
Next, using \eqref{E-39} and \eqref{E-41}, the above equality simplifies to
\begin{align} \label{E-41a}
& 2(\nabla_X {F})(Z,Y)+ 2(\nabla_{f X}{F})(Z,fY)  
 = 4\,\overline{\eta}(Z) g(QX,QY) 
 -\sum\nolimits_j\big[2\,\eta^j(Y)g((Q+h_j)X,QZ) \notag \\ 
&+4\,\eta^j(X)\eta^j(Y)\overline{\eta}(Z)
 -2\,\overline{\eta}(Y)\eta^j(X)\eta^j(Z) 
 -\eta^j(Y)g((Q+h_j)Z, \widetilde Q X)\notag 
 +\eta^j(Z)g((Q+h_j)Y, \widetilde Q X)\\ 
&-\eta^j(X)g((Q+h_j)Z, \widetilde Q Y)\big] 
 +(\nabla_X {F})(Y,\widetilde QZ)-(\nabla_Y {F})(Z,\widetilde QX)-(\nabla_Z {F})(X,\widetilde QY).
\end{align}
Observe that using \eqref{E-nS-10}, and then the anti-symmetry of $f$ and $\nabla Q=\nabla \widetilde Q$, we have 
\begin{align} \label{E-Phi}
 (\nabla_X {F})(Y,\widetilde QZ)= (\nabla_X {F})(\widetilde QY,Z){+}\sum\nolimits_i\big[
 \eta^j(Z)g((Q+h_j)X,\widetilde QY)
 +\eta^j(Y)g((Q+h_j)X,\widetilde QZ)\big] 
\end{align}
for all $X,Y,Z \in\mathfrak{X}_M$.
Therefore, from the above and $d{F}=0$, the equation \eqref{E-41a} simplifies as
\begin{align}\label{E-49new}
\notag
 (\nabla_X {F})(Z,Y) & +(\nabla_{f X} {F})(Z,f Y) 
 = 2\,\overline{\eta}(Z) g(QX,QY) 
 -\sum\nolimits_j\big[\eta^j(Y)g((Q+h_j)X,QZ) \\
\notag
 & +2\,\eta^j(X)\eta^j(Y)\overline{\eta}(Z) 
 -\eta^j(X)\overline{\eta}(Y)\eta^j(Z) +\tfrac12\,\eta^j(Z)\,g((Q+h_j)Y, \widetilde Q X)\big] \\
 &+\tfrac12\big[(\nabla_X {F})(Y,\widetilde QZ) 
 +(\nabla_{\widetilde Q X} {F})(Y,Z) \big] ,
\end{align}
that is equivalent to \eqref{E-04}.
\end{proof}

The subsequent result generalizes \cite[Lemma~2.1]{DTL-2006} and is used in Theorem~\ref{P-07} and 
Proposition~\ref{prop-5.1}.

\begin{proposition}\label{prop-5.2}
The curvature tensor of a w.a.$\,\mathcal{S}$-manifold with \eqref{Eq-2assump} 
and \eqref{E-nS-10} satisfies
\begin{align}\label{E-05b}
 & g(R_{\xi_i,QX}Y, Z) -g(R_{\xi_i,X}f Y, f Z) 
 +g(R_{\xi_i,f X}Y, f Z) -g(R_{\xi_i,f X}f Y, Z) 
 = 2\,(\nabla_{\widetilde h_iX} {F})(Y,Z) \notag \\ 
 & +2\,\overline{\eta}(Z)g((Q+Q\widetilde h_i)X,QY)
 {-}2\,\overline{\eta}(Y)g((Q+Q\widetilde h_i)X, QZ) \notag 
 {-}\sum\nolimits_j\eta^j(X) \big[2\,\eta^j(Y) \overline{\eta}(Z) {-} 2\,\overline{\eta}(Y)\eta^j(Z)\\ 
 & -g([\widetilde h_i,\widetilde h_j]Z,\widetilde Q^2QY)\big] 
 +\tfrac12\big[3(\nabla_{\widetilde Q\widetilde h_iX} {F})(Y,Z) 
 -(\nabla_{\widetilde QY} {F})(Z,\widetilde h_iX) 
 +(\nabla_Y {F})(\widetilde h_iZ, \widetilde QX)  
  \notag \\ 
 & -(\nabla_{\widetilde QY} {F})(\widetilde h_i Z, X) 
 -(\nabla_{\widetilde QZ} {F})(\widetilde h_iX,Y)
 +(\nabla_Z {F})(\widetilde QX,\widetilde h_iY) 
 -(\nabla_{\widetilde QZ} {F})(X,\widetilde h_iY)\big].
\end{align}
\end{proposition}

\begin{proof}
Using \eqref{E-xi-n}, proved using the condition \eqref{Eq-2assump}, the curvature tensor $R_{Y,Z}\,\xi_i$ is given by
\begin{align}\label{E-42}
 R_{Y,Z}\, \xi_i =-(\nabla_Y f )Z+(\nabla_Z f )Y
 -(\nabla_Y f\widetilde h_i)Z 
 +(\nabla_Z f\widetilde h_i)Y,
\end{align}
Taking the scalar product of \eqref{E-42} with $X$ and 
using \eqref{dPhi} (since $d{F}=0$), we get 
\begin{align}\label{E-05}
 g(R_{\xi_i,X}Y ,Z)=-(\nabla_X {F})(Y,Z)-g(X,(\nabla_Y f \widetilde h_i)Z)+ g(X, (\nabla_Z  f  \widetilde h_i)Y).
\end{align}
Using \eqref{E-05}, we have 
\begin{align}\label{E-42-tilde}
 g(R_{\xi_i,X}Y, Z) &-g(R_{\xi_i, X}f Y, f Z) 
 + g(R_{\xi_i,f X}Y, f Z)
 +g(R_{\xi_i,f X}f Y, Z) \notag \\
 &= A(X,Y,Z) + \widetilde B_i(X,Y,Z) - \widetilde B_i(X,Z,Y) ,
\end{align}
where the operators $A$ and $\widetilde B_i\ (1\le i\le s)$ are defined as follows:
\begin{align*}
 & A(X,Y,Z) = -(\nabla_X {F})(Y,Z) +(\nabla_X {F})(f Y,f Z) -(\nabla_{f X} {F})(Y, f Z)  -(\nabla_{f X} {F})(f Y,Z), \\
 & \widetilde B_i(X,Y,Z) = 
 g(X, (\nabla_{f Y}f\widetilde h_i)f Z)
 -g(X, (\nabla_Yf\widetilde h_i)Z) 
 - g(f X, (\nabla_Yf\widetilde h_i)f Z) 
 - g(f X, (\nabla_{f Y}f\widetilde h_i)Z).
\end{align*}
Using 
 \eqref{E-39}, 
 \eqref{E-41} 
 and 
 \eqref{E-49new}
of Proposition~\ref{P-05}, 
and Proposition~\ref{L-cond-equiv},
the operator $A$ reads as
\begin{align}\label{eq-a1}
\notag
A(X,Y,Z) =& \,\big\{
 (\nabla_X {F})(Y,Z)+(\nabla_{X} {F})(fY,fZ)\big\}
+\big\{(\nabla_{fX} {F})(fY,Z)
-(\nabla_{fX} {F})(Y, fZ)\big\} \\
\notag
 =&\,2\,\overline{\eta}(Z) g(QX, QY) -2\,\overline\eta(Y) g(QX, QZ)  
 +2\sum\nolimits_j
 \eta^j(X)\big\{\overline\eta(Y)\eta^j(Z) 
 -\overline{\eta}(Z)\eta^j(Y)\big\} \\
 & +(\nabla_{\widetilde Q X} {F})(Y,Z) ,
\end{align}
and the operator $\widetilde B_i$ reads as
\begin{align*}
 &\widetilde B_i(X,Y,Z) 
 =\, g( X, f(\nabla_Y \widetilde Q \widetilde h_i)Z) 
 -g( X,(\nabla_{Y} f) \widetilde h_iZ) 
+g(X, (\nabla_{f Y}Q \widetilde h_i )Z)  +g(X, \widetilde h_i f (\nabla_{f Y}f) Z) \\ 
 & +g(Q X, \widetilde h_i(\nabla_Y f) Z)+g( X, f(\nabla_{f Y}f) \widetilde h_i Z)
 -g(Q X, (\nabla_{f Y} \widetilde h_i) Z)
 +\sum\nolimits_j\eta^j(X)\eta^j((\nabla_{f Y} \widetilde h_i) Z).
\end{align*} 
Using \eqref{E-nS-10}, we get 
\begin{align*}
g(f X, (\nabla_Y \widetilde Q\,\widetilde h_i)Z) &= 
g(f X, \widetilde Q( \nabla_Y \widetilde h_i)Z), \\ 
 g(X, (\nabla_{f Y}Q\,\widetilde h_i) Z) &= 
 g(X,  Q(\nabla_{f Y} \widetilde h_i )Z) 
 -\sum\nolimits_j 
 \eta^j(X)g((Q-Q \widetilde h_j) Y,\widetilde Q Q \widetilde h_iZ).
\end{align*} 
These two equations allow us to simplify 
$\widetilde B_i(X,Y,Z)$ as follows:
\begin{align*}
 \widetilde B_i(X,Y,Z) 
 =& -g(X,(\nabla_{Y}f)\widetilde h_iZ) 
 +g(Q X, \widetilde h_i(\nabla_Y f) Z) 
 +g(X, \widetilde h_i f (\nabla_{f Y}f) Z) \\ 
 &\,
 +g(X, f(\nabla_{f Y}f)\widetilde h_i Z) 
 +\sum\nolimits_j\eta^j(X)\eta^j((\nabla_{fY} \widetilde h_i) Z) \\
 & -g(\widetilde Q f X, (\nabla_Y \widetilde h_i)Z) 
  -\sum\nolimits_j 
 \eta^j(X)g((Q-Q\widetilde h_j)Y,\widetilde Q Q \widetilde h_iZ) .
\end{align*}
To further simplify $\widetilde B(X,Y,Z)$, we compute $f(\nabla_{f Y} f)Z$ using 
Definition~\ref{Def-01}, \eqref{E-xi-n}, \eqref{E-nS-10b} and \eqref{E-04}, 
\begin{align} \label{e-phi}
 f (\nabla_{fY} f)Z
 & =(\nabla_{f Y}f ^2)Z -(\nabla_{f Y} f)f Z\notag \\
 & = -\big\{ g(QY,QZ)-\sum\nolimits_j\eta^{j}(Y)\eta^j(Z)\big\}
 \overline{\xi} +2\,\overline{\eta}(Z)\big\{Q^2Y
 -\sum\nolimits_{j}\eta^j(Y)\xi_j\big\}\notag \\
 & \quad -\sum\nolimits_{j}g(Q \widetilde h_jY,QZ)\xi_j 
 +(\nabla_Y f)Z +P(Y,Z) 
 +\tfrac{1}{2}\sum\nolimits_{j}
  g((Q+Q \widetilde h_j)Z,\widetilde QY)\xi_j ,
\end{align}
where  
$P(X,Y)=\tfrac{1}{2}\big\{\widetilde Q(\nabla_X f) Y 
 +(\nabla_{\widetilde QX} f) Y \big\}$.
Using the equality $\eta^i((\nabla_{f Y} \widetilde h_j)Z)
= g(Q\widetilde h_jZ, \widetilde h_iY-Y)$ for 
$i,j=1,\dots, s$, and \eqref{e-phi}, 
we simplify $\widetilde B_i(X,Y,Z)$ to 
\begin{align}\label{Eq-b1}
\notag
 &\widetilde B_i(X,Y,Z) = 2\,g(\widetilde h_iX, (\nabla_{Y}f) Z) 
 +2\,\overline{\eta}(Z)g( Q \widetilde h_iX, QY) -2\,\overline{\eta}(X)g(QY, Q \widetilde h_iZ)+g(\widetilde{Q}X, f (\nabla_Y \widetilde h_i)Z) \\ & \quad
 +g(\widetilde{Q}X, \widetilde h_i(\nabla_{Y}f) Z)  {+}g(P(Y,Z),\widetilde h_iX){+}g(P(Y, \widetilde h_iZ),X) 
 +\tfrac{1}{2}\sum\nolimits_j\eta^j(X) g((Q+Q \widetilde h_j)Y,\widetilde QZ)  .
\end{align} 
Finally, we simplify \eqref{E-42-tilde} using $d{F}=0$, \eqref{E-Phi}, \eqref{eq-a1}, and 
\eqref{Eq-b1} as follows:
\begin{align}\label{E-ABB} 
 & A(X,Y,Z) + \widetilde B_i(X,Y,Z) - \widetilde B_i(X,Z,Y) \notag\\ 
 & \quad = 2\,(\nabla_{\widetilde h_iX} {F})(Y,Z) 
 +2\,\overline{\eta}(Z)g((Q+Q \widetilde h_i)X,QY)- 2\,\overline{\eta}(Y)g((Q+Q \widetilde h_i)X, QZ) \notag \\ 
 & \quad -2\,\sum\nolimits_j\eta^j(X) \big[\eta^j(Y) \overline{\eta}(Z)- \overline{\eta}(Y) \eta^j(Z) \big] + (\nabla_{\widetilde QX} {F})(Y,Z) \notag   - (\nabla_{\widetilde Q \widetilde h_iX} {F})(Y,Z) \notag \\  
 & \quad  + g(\widetilde QX, f (\nabla_Y \widetilde h_i)Z-f (\nabla_Z \widetilde h_i)Y)\notag + g(P(Y,Z),\widetilde h_iX)+g(P(Y, \widetilde h_iZ),X) - g(P(Z,Y),\widetilde h_iX) \\ 
 & \quad-g(P(Z, \widetilde h_iY),X) 
 +\sum\nolimits_j\eta^j(X)
 g([\widetilde h_i,\widetilde h_j]Z, \widetilde Q^2QY).
\end{align}
Using  $d{F}=0$ and \eqref{E-Phi}, we simplify the terms of \eqref{E-ABB} containing the tensor $P$ and get the following:
\begin{align*}
&2\,g(P(Y,Z), \,\widetilde h_iX)  +2\,g(P(Y,\widetilde h_iZ),X)  -2\,g(P(Z,Y), \widetilde h_iX) -2\,g(P(Z,\widetilde h_iY),X) \\
 & =(\nabla_{\widetilde Q\widetilde h_i X}{F})(Y,Z) 
 -(\nabla_{\widetilde QY} {F})(Z,\widetilde h_iX)  
 + (\nabla_Y {F})(\widetilde QX,\widetilde h_iZ) 
 -(\nabla_{\widetilde QY} {F})(\widetilde h_iZ,X) \\
 &\quad  +(\nabla_{\widetilde QZ} {F})(Y,\widetilde h_iX)
 -(\nabla_Z {F})(\widetilde QX,\widetilde h_iY) +(\nabla_{\widetilde QZ}{F})(\widetilde h_iY,X).
\end{align*}
From \eqref{E-05} we have 
\begin{align*}
 g(\widetilde Q X, f(\nabla_Z\widetilde h_i)Y
 {-}f(\nabla_Y\widetilde h_i)Z)   
 = g(\widetilde Q X, (\nabla_Z\,f\widetilde h_i)Y 
 {-}(\nabla_Y\,f\widetilde h_i)Z)
 {-}g(\widetilde Q X, (\nabla_Z\,f)\widetilde h_i Y 
 {-}(\nabla_Y\,f)\widetilde h_i Z) \\ 
 \quad\quad = g(R_{\xi_i,\widetilde Q X}Y, Z)
 + (\nabla_{\widetilde QX} {F})(Y,Z)
 -(\nabla_Z{F})(\widetilde Q X,\widetilde h_i Y) 
 +(\nabla_Y{F})(\widetilde Q X,\widetilde h_i Z).
 \end{align*}
Using two expressions above, we simplify \eqref{E-ABB} to the following:
\begin{align} \label{E-CURV}
  &A(X,Y,Z) + \widetilde B_i(X,Y,Z) - \widetilde B_i(X,Z,Y)\notag \\  
 \notag
 &  = 2\,(\nabla_{\widetilde h_iX} {F})(Y,Z)  +2\,\overline{\eta}(Z)g((Q+Q \widetilde h_i)X,QY)- 2\,\overline{\eta}(Y)g((Q+Q \widetilde h_i)X, QZ) \notag  \\ 
 & \quad -\sum\nolimits_j\eta^j(X) \big[2\,\eta^j(Y) \overline{\eta}(Z)- 2\,\overline{\eta}(Y) \eta^j(Z) -g([\widetilde h_i,\widetilde h_j]Z, \widetilde Q^2QY)\big]  
 -g(R_{\xi_i,\widetilde Q X}Y, Z)  \notag
 \\ \notag 
 &\quad+ \tfrac{1}{2} \big\{3(\nabla_{\widetilde Q\widetilde h_iX} {F})(Y,Z) -(\nabla_{\widetilde QY} {F})(Z,\widetilde h_iX) 
 +(\nabla_Y {F})(\widetilde h_iZ, \widetilde QX)  -(\nabla_{\widetilde QY} {F})(\widetilde h_i Z, X)
    \notag \\ 
 & \quad  -(\nabla_{\widetilde QZ} {F})(\widetilde h_iX,Y)
 +(\nabla_Z {F})(\widetilde QX,\widetilde h_iY) 
 -(\nabla_{\widetilde QZ} {F})(X,\widetilde h_iY)\big\}.
\end{align}
The required equality \eqref{E-05b} follows from \eqref{E-42-tilde} and \eqref{E-CURV}.
\end{proof}

The following statement plays a key role in the proof of Theorem~\ref{T-04} and in Example~\ref{C-K-S2}.

\begin{proposition}\label{prop-5.1}
Let a w.a.$\,\mathcal{S}$-manifold 
\((M^{2n+s},f,Q,\xi_i,\eta^i,g)\) with the $f$-$(\kappa,\mu)$-nullity property \eqref{E-k-mu-nullity} satisfy the conditions  \eqref{Eq-2assump} and \eqref{E-nS-10}. 
If $\kappa<1$, then
$\widetilde h_1=\ldots=\widetilde h_s$, and we have the following:
\begin{align}
 (\nabla_X {F})(Y,Z) & = 0  
 \qquad (X,Z\in {\cal D_-}, \ Y\in {\cal D_+}),
 \label{Eq-phi03} \\
 (\nabla_X {F})(Y,Z) &=0 \qquad (X,Y \in {\cal D}_{+},\ Z \in {\cal D} ),
 \label{Eq-phi02}\\
 \label{Eq-phi001}
 (\nabla_X {F})(Y,Z) &=0 \qquad 
 (X,Y,Z\in {\cal D}_-),
\end{align}
where ${\cal D}_\pm$ are $n$-dimensional eigen-distributions of the eigenvalues $\pm\sqrt{1-\kappa}$
of~$\widetilde h_i$.
\end{proposition}

\begin{proof}
By \eqref{E-k-mu-nullity} and \eqref{E-37}, we obtain 
$h_i^2=(\kappa-1)\,Qf^2$. 
Since $\widetilde h_i$ is self-adjoint, we have $\widetilde h_1=\ldots=\widetilde h_s$.
Next, from \eqref{E-05b} of Proposition~\ref{prop-5.2}, by considering $X,Y,Z \in {\cal D}$ and applying~\eqref{E-RXY-xi}, we obtain
\begin{align}\label{Eq-phi01}
  4(\nabla_{\widetilde hX} {F})(Y,Z) 
 &= (\nabla_{\widetilde QY} {F})(Z,\widetilde hX) 
 +(\nabla_{\widetilde QY}{F})(\widetilde h Z,X)  
 +(\nabla_{\widetilde QZ} {F})(\widetilde hX,Y)  
 +(\nabla_{\widetilde QZ} {F})(X,\widetilde hY)
 \notag \\ 
 &-3(\nabla_{\widetilde Q\widetilde hX} {F})(Y,Z) 
 -(\nabla_Y {F})(\widetilde hZ, \widetilde QX)  
 -(\nabla_Z {F})(\widetilde QX,\widetilde hY) \quad (X,Y,Z \in {\cal D}) ,   
\end{align}
where we set $\widetilde h=\widetilde h_i$.
Taking $X,Z\in {\cal D^-}$ and $Y\in {\cal D^+}$ in \eqref{Eq-phi01}, and using $d{F}=0$, we~find 
\begin{align} \label{Eq-N}
 4\, (\nabla_X {F})(Y,Z) 
 &= -3(\nabla_{\widetilde Q X} {F})(Y,Z)
 +2\,(\nabla_{\widetilde Q Y} {F})(Z,X)
 -(\nabla_{Y} {F})(\widetilde Q Z,X) 
 +(\nabla_{Z} {F})(\widetilde Q X,Y) \notag \\
 & =-2\,\big[(\nabla_{\widetilde Q X} {F})(Y,Z) 
 +(\nabla_{X}{F}) (\widetilde QY,Z)\big].
\end{align} 
Cyclically changing $X,Y,Z$ and using the above in conjunction with $d{F}=0$, we acquire
\begin{align*}
 0 &= -2\,d{F}(X,Y,Z)= - 2\, [(\nabla_X {F})(Y,Z)
 +(\nabla_Y {F})(Z,X)+ (\nabla_Z {F})(X,Y)] \\
 & = (\nabla_{\widetilde Q X} {F})(Y,Z) 
 + (\nabla_{X}{F}) (\widetilde QY,Z) 
 + (\nabla_{\widetilde Q Y} {F})(Z,X) 
 + (\nabla_{Y}{F}) (\widetilde QZ,X) \\ 
 & \quad + (\nabla_{\widetilde Q Z} {F})(X,Y) 
 + (\nabla_{Z}{F}) (\widetilde QX,Y) 
 = (\nabla_{\widetilde Q X} {F})(Y,Z) 
 - (\nabla_{X}{F}) (\widetilde QY,Z).
\end{align*}
Using this in \eqref{Eq-N}, gives $(\nabla_X\,{F})(QY,Z)=0$, which implies \eqref{Eq-phi03} as $Q$ is a bijection on ${\cal D}_+$. 
Taking $X,Y,Z\in {\mathcal D}^+ $(or $X,Y,Z\in {\cal D}_-$)  in \eqref{Eq-phi01}, we find 
\begin{align*}
 4(\nabla_{X}{F})(Y,Z) 
 &= -3(\nabla_{\widetilde QX}{F})(Y,Z) 
 +2(\nabla_{\widetilde QY}{F})(Z,X) +2(\nabla_{\widetilde QZ}{F})(X,Y) \\ 
 &\quad-(\nabla_{Z}{F})(\widetilde Q X,Y)
 -(\nabla_{Y}{F})(Z,\widetilde Q X). 
\end{align*}
Using this and \eqref{dPhi}, we get $(\nabla_{X}{F})(QY,Z)=0$, which indicates that $(\nabla_{X}{F})(Y,Z)=0$ for $X,Y,Z\in {\mathcal D}^+$(or $X,Y,Z\in {\cal D}_-$), that is, \eqref{Eq-phi02} is true. The proof of \eqref{Eq-phi001}  is similar.
\end{proof}

\section{Conclusion}
\label{sec:concl}


In this work, we investigated curvature conditions on w.a.\,$\cal S$-manifolds, focusing on the $f$-$(\kappa, \mu)$-nullity condition and its special case $R_{X,Y}\,\xi=0$, and established 
results that broaden the classical~theory.
We~first showed that the curvature identity characterizing $\cal S$-manifolds remains invariant under the partial Ricci flow (PRF) when the weak structure satisfies mild compatibility conditions. 
This allowed us to obtain a {dynamical characterization of $\cal S$-structure}: starting from a w.a.\,$\cal S$-structure satisfying the curvature condition of $\cal S$-manifolds, the PRF evolves the structure exponentially fast toward a genuine $\cal S$-structure. This extends the dynamical results previously known for $f$-K-contact manifolds and shows that the PRF is a natural tool for detecting and stabilizing $\cal S$-structures.
Next,~we analyzed w.a.\,$\cal S$-manifolds satisfying the $f$-$(\kappa, \mu)$-nullity condition. 
For $\kappa=1$, we proved that such manifolds admit a $\cal D$-adapted deformation to an $\cal S$-structure, again via the PRF. For $\kappa<1$, we established the existence of {bi-Legendrian structures}: the contact distribution decomposes into conjugate Legendrian eigen-distributions of the structural tensor $\widetilde h_i$, each defining a totally geodesic foliation. 
This gene\-ralizes results of Cappelletti Montano and Di Terlizzi to the weak metric setting and demonstrates that w.a.\,$\cal S$-geometry naturally accommodates rich 
bi-Legendrian behavior.
We established a {sca\-ling correspondence} between almost $\cal S$-structures and w.a.\,$\cal S$-structures satisfying the $f$-$(\kappa,\mu)$-nullity condition.
Finally, in the case $\kappa=\mu=0$, we proved a {splitting theorem} in which one factor is flat
and the another has positive constant sectional curvature. 
This extends Di Terlizzi’s splitting results for almost $\cal S$-manifolds and shows that the weak metric framework retains strong geometric rigidity under vanishing curvature conditions.
Our results demonstrate that w.a.\,$\cal S$-manifolds preserve the essential geometric properties of classical $\cal S$-geometry while offering a broader and more adaptable setting. The~PRF emerges as a powerful tool for understanding the evolution, rigidity, and classification of these manifolds. 
We~expect that the techniques developed here will be useful in further studies of weak metric structures, the dynamics of Legendrian and contact foliations, weak $f$-$(\kappa,\mu)$-manifolds in CR, complex, and Lorentzian settings, non-symmetric metric affine geometry, 
having 
applications in mathematical physics (see~\cite{BC-26,RZ-1}).


\section*{\vskip-10mm\noindent Declarations}

\begin{itemize}
\item\vskip-2.1mm Funding:
The first author would like to thank the University Grant Commission (UGC) (Grant number 211610029330) for providing financial support.

\item\vskip-2.1mm {Conflict of interest:} The authors declare no conflicts of interest.

\item\vskip-2.1mm {Ethics approval and consent to participate:} The work is original, not under consideration elsewhere, and approved by all authors.

\item\vskip-2.1mm {Author contribution:} All authors contributed equally to this work.

\item\vskip-2.1mm {Data availability statement:}
The manuscript has no associated data.
\end{itemize}


\end{document}